\documentclass[11pt]{article}
\usepackage{amssymb,amsmath,amsfonts,amsthm}
\usepackage{graphicx}
\usepackage{subfigure}
\usepackage{makeidx}
\usepackage{multicol}
\usepackage{tikz}
\usetikzlibrary{positioning}
\usepackage[top=1in, bottom=1.25in, left=1.25in, right=1.25in]{geometry}

\newtheorem{theorem}{Theorem}
\newtheorem{example}{Example}
\newtheorem{definition}[theorem]{Definition}

\newtheorem{proposition}[theorem]{Proposition}

\newcommand{\E}{\mathbb{E}}

\newcommand{\cH}{\mathcal{H}}
\newcommand{\cI}{\mathcal{I}}

\renewcommand{\P}{\mathbb{P}}

\newcommand{\R}{\mathbb{R}}

\newcommand{\cX}{\mathcal{X}}

\newcommand{\cY}{\mathcal{Y}}
\newcommand{\indep}{\mbox{\,$\perp\!\!\!\perp$\,}}

\begin{document}
\author{Piotr Zwiernik}
\title{\textbf{Latent tree models}}
\date{}
\maketitle

\begin{abstract}
	\textit{Latent tree models} are graphical models defined on trees, in which only a subset of variables is observed. They were first discussed by Judea Pearl as tree-decomposable distributions to generalise star-decomposable distributions such as the \emph{latent class model}. Latent tree models, or their submodels, are widely used in: phylogenetic analysis, network tomography, computer vision, causal modeling, and data clustering. They also contain other well-known classes of models like hidden Markov models, Brownian motion tree model, the Ising model on a tree, and many popular models used in phylogenetics. We offer here a concise introduction to the theory of latent tree models. We emphasise the role of \emph{tree metrics} in the structural description of this model class, in designing learning algorithms, and in understanding fundamental limits of what and when can be learned. 

\end{abstract}

\tableofcontents

%\section{Star-decomposable distributions}

\section{Basics}

In this section we define latent tree models and provide motivation to work with this model class. We present Gaussian and general Markov models as subclasses of latent tree models that admits tractable and rigorous analysis.

\subsection{Definitions}

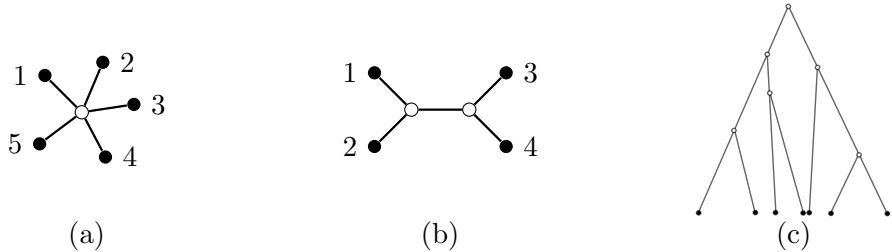
\begin{figure}[htp!]
	\begin{minipage}{.3\textwidth}\centering
	\tikzstyle{vertex}=[circle,fill=black,minimum size=5pt,inner sep=0pt]
\tikzstyle{hidden}=[circle,draw,minimum size=5pt,inner sep=0pt]
\tikzstyle{vertex}=[circle,fill=black,minimum size=5pt,inner sep=0pt]
\tikzstyle{hidden}=[circle,draw,minimum size=5pt,inner sep=0pt]
  \begin{tikzpicture}[scale=.7]
  \node[vertex] (1) at (-.7,.7)  [label=left:$1$] {};
    \node[vertex] (2) at (.4,.95) [label=right:$2$]{};
    \node[vertex] (3) at (.99,.15) [label=right:$3$]{};
    \node[vertex] (4) at (.45,-.85) [label=right:$4$]{};
    \node[vertex] (5) at (-.8,-.6)  [label=left:$5$] {};
    \node[hidden] (a) at (0,0) {};
\coordinate (c) at (.8,-0.1);
\coordinate (d) at (-.5,-0.5);
    \draw[line width=.3mm] (a) to (1);
    \draw[line width=.3mm] (a) to (2);
    \draw[line width=.3mm] (a) to (3);
        \draw[line width=.3mm] (a) to (4);
        \draw[line width=.3mm] (a) to (5);
  \end{tikzpicture}		\end{minipage}
\begin{minipage}{.3\textwidth}\centering
\tikzstyle{vertex}=[circle,fill=black,minimum size=5pt,inner sep=0pt]
\tikzstyle{hidden}=[circle,draw,minimum size=5pt,inner sep=0pt]
\begin{tikzpicture}[scale=.7]
  \node[vertex] (1) at (-.7,.7)  [label=left:$1$] {};
    \node[vertex] (2) at (-.7,-.7) [label=left:$2$]{};
    \node[vertex] (3) at (1.8,.7) [label=right:$3$]{};
    \node[vertex] (4) at (1.8,-.7) [label=right:$4$]{};
    \node[hidden] (a) at (0,0) {};
    \node[hidden] (b) at (1.1,0) {};
    \draw[line width=.3mm] (a) to (1);
    \draw[line width=.3mm] (a) to (2);
    \draw[line width=.3mm] (b) to (3);
        \draw[line width=.3mm] (b) to (4);
        \draw[line width=.3mm] (a) to (b);
  \end{tikzpicture}  \end{minipage}
  \begin{minipage}{.3\textwidth}\centering
\includegraphics[scale=.15]{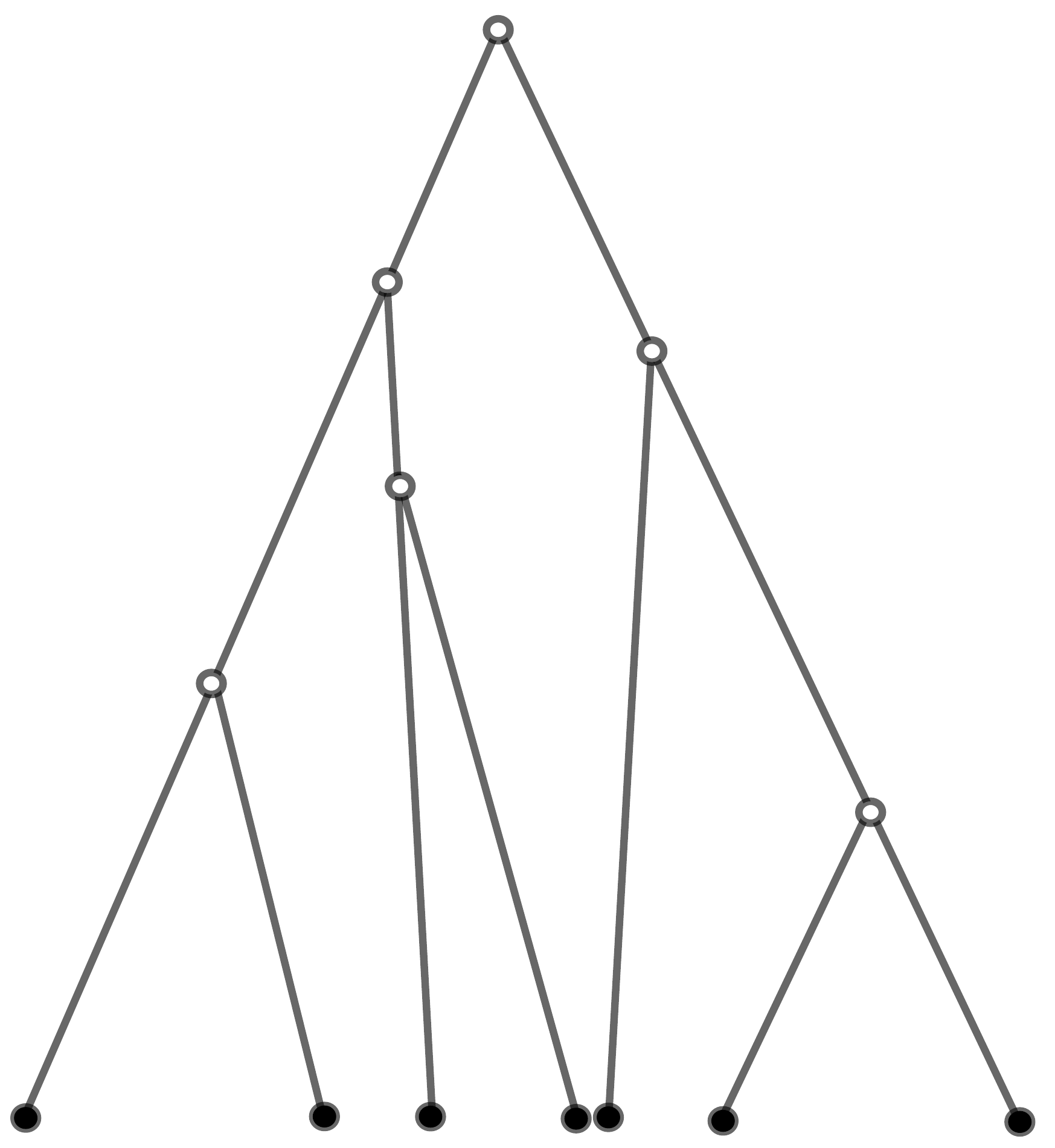}
  \end{minipage}\\
    \begin{minipage}{.3\textwidth}\centering
  (a)
  \end{minipage}  \begin{minipage}{.3\textwidth}\centering
  (b)
  \end{minipage}  \begin{minipage}{.3\textwidth}\centering
  (c)
  \end{minipage}
  \caption{(a) An undirected star tree with five leaves, (b) a quartet tree, (c) a binary rooted tree.}\label{fig:6}
\end{figure}

A \emph{tree} is an undirected graph without cycles. A \emph{leaf} of $T$ is a vertex of degree one, an \emph{internal vertex} is a vertex which is not a leaf, and an \emph{inner edge} is an edge whose both ends are internal vertices. Given a tree $T$ define a \emph{rooted tree} as a directed graph obtained from $T$ by picking one of its vertices $r$ and directing all edges away from $r$. The vertex $r$ is called the \emph{root}. Trees will be always leaf-labeled with the labelling set $\{1,\ldots,m\}$, where $m$ is the number of leaves. An undirected tree is \emph{trivalent} if each internal vertex has degree precisely three. A rooted tree is a \emph{binary rooted tree} if each internal vertex has precisely two children. In many applications rooted trees are depicted without using arrows, where direction is made implicit by drawing the root on the top and the leaves on the bottom; see Figure~\ref{fig:6}(c). Two special types of undirected trees are: a \emph{star tree} with one internal vertex and a trivalent tree on four leaves called a \emph{quartet tree}; see Figure~\ref{fig:6}(a) and (b). A \emph{forest} is a collection of trees. Forests here are also leaf-labeled with  the labelling set is $\{1,\ldots,m\}$, which means that each tree in this collection is leaf-labeled and the corresponding collection of labelling sets forms a set partition of $\{1,\ldots,m\}$. We define three graph operations on trees (forests). \emph{Removing an edge} means removing that edge from the edge set. \emph{Contracting an edge} $u-v$ means removing $u,v$ from the vertex set, adding a new vertex $w$ and edges such that $w$ is adjacent to all vertices which were adjacent to $u$ or $v$. \emph{Suppressing a vertex of degree two} means removing that vertex and replacing the two edges incident to that vertex by a single edge.

A latent tree model is defined as follows. Let $Y=(Y_v)_{v\in V}$ be a random vector with coordinates indexed by the vertices $V$ of a tree $T$ and with values in the state space \mbox{$\cY=\prod_{v\in V}\cY_v$}. Suppose that the density of $Y$, denoted by $p(y)$, lies in the graphical model over $T$. This means that $p$ factorises acording to the tree $$p(y)\;\;=\;\;\prod_{u-v} \Psi_{uv}(y_u,y_v),$$ 
where for each edge $u-v$ of $T$ the function $\Psi_{uv}$ is a nonnegative function. We call such a model a \emph{fully-observed} tree model. Denote the set of leaves by $W\subset V$ and let $m:=|W|$. Write $X:=Y_W$, $\cX:=\cY_W$, $H:=Y_{V\setminus W}$, and $\cH:=\cY_{V\setminus W}$. The latent tree model $M=M(T,\cY)$ is the family of marginal distributions of $p(y)$ over the leaves $W$:
\begin{equation}\label{eq:margin}
p_W(x)\;\;:=\;\;\int_{\cH} p(x,h) \,{\rm d}h.	
\end{equation}
In other words, the internal vertices represent unobserved random variables. The above definition extends to situations where some internal vertices are also observed. However, these seemingly more general situation does not lead to any new family of distributions; for Gaussian and general Markov models this is shown in Theorem \ref{th:binaryred}. 
 
% Similarly, we can consider more generally models over forests, where each tree component is treated as independent. Although we do not discuss it here, these models also are naturally submodels of models over trees\textcolor{blue}{CITATION}. 
%  
Consider now a Bayesian network on a rooted tree obtained from $T$. This is a model of distributions that factorise as follows
\begin{equation}\label{eq:fullyfactor}
p(y)\;\;=\;\;p_{r}(y_{r})\prod_{u\to v}p_{v|u}(y_{v}|y_{u}) \quad\mbox{for all }y\in \cY	.
\end{equation}
By standard results on Markov equivalence of directed acyclic graphs \cite{frydenberg1990,vermapearl91}, this  model coincides with the latent tree model on $T$ for any choice of the root location. The fully-observed tree model is then fully characterized by the root distribution $p_r(y_r)$ and by the conditional distributions $p_{v|u}(y_v|y_u)$ for all edges $u\to v$ in the rooted tree. The parameterisation of the corresponding latent tree model is induced by taking the margin over $\cX$ as in~(\ref{eq:margin}). 

Although the definition of latent tree models is fairly general, here we make additional assumptions on the state space $\cY$ and possible distributions. If $\cY$ is a finite set then we call the corresponding latent tree model \emph{discrete}. \emph{Gaussian latent tree models} are latent tree models for which the vector $Y$ is jointly Gaussian. In the Gaussian case we typically assume that the mean of $Y$ is known and equal to zero. Mixed cases when $\cH$ is finite and $X$ is Gaussian conditionally on $H$ are also popular.  

The parameters of latent tree models are of two kinds. The underlying tree $T$ is the \emph{discrete parameter}. The \emph{continuous parameter} $\theta$ is the parameter specifying the root distribution and the conditional distributions for each edge. To make the parameters explicit, we write $p(y)=p(y;T,\theta)$ and $p_W(x)=p_W(x;T,\theta)$. Formally, the state space $\cH$ of the unobserved part of the vector $Y$ is also a parameter of a latent tree model. In certain applications finding $\cH$ can be important; see~\cite[Section 3.2.4]{mourad2013survey}. Here we always assume that $\cH$ is fixed. With this convention, there are three main learning tasks related to latent tree models:
\begin{enumerate}
	\item[(L1)] Given a sample of size $n$ from a latent tree model estimate the underlying tree $T$.
	\item[(L2)] Given a sample of size $n$ from a latent tree model on a fixed tree $T$ estimate the continuous parameter $\theta$.
	\item[(L3)] Given a fully specified latent tree model and a single sample at the observed vertices, infer the states at the unobserved vertices.
\end{enumerate}
There are several fundamental questions related to the first two learning tasks that we are going to address in this exposition. We will not discuss here the learning task (L3). If the model is completely specified then, by~(\ref{eq:fullyfactor}), we have access to the full distribution over $\cY$. In this case the learning task (L3) reduces to the sum-product algorithm discussed, for example, in  \cite[Section 8.4]{bishop2006pattern}. 

%Although the analysis of this model class is complicated, we will show that as long as correlations between variables are high both learning task can be achieved with a bounded error. Te complication occurs when some correlations are low, in which case some phase transition occurs and suddenly standard statistical theory does not apply. We will discuss these issues in more detail in Section~\ref{sec:theory}.

\subsection{Motivation and  applications}\label{sec:apps}

Latent tree models form the most tractable family of Bayesian networks with unobserved variables, which can be used to model dependence structures when unobserved confounders are expected; see, for example,~\cite[Section 2]{pearl2000} or~\cite{evans2015graphs}. However, there are several other reasons why latent tree models become popular across sciences. We distinguish three main types of applications. 

First, latent tree models represent a larger family of probability distributions than fully-observed tree models but retain some of their computational advantages, which is particularly important in high-dimensional settings. From~(\ref{eq:fullyfactor}) it is clear that having an estimator $\hat\theta$ of model parameters we obtain an estimator of the fully observed distribution $p(y;\hat\theta)$ and so we can very efficiently compute various marginal distributions in the model using the sum-product algorithm; see, for example, \cite[Section 8.4]{bishop2006pattern}. Using the max-product algorithm it is also possible to efficiently infer the unobserved states. 

Second, a rooted tree can represent evolutionary processes with the root representing the common ancestor and the leaves representing extant species. This makes latent tree models useful in phylogenetic analysis~\cite{felsenstein2004inferring, semple2003pol}. In this context, the data typically consist of $m$ aligned DNA sequences of length $n$, where each site in the sequence is treated as an independent realisation of the vector $X$. In this case all state spaces are of the form $\{A,C,T,G\}$ or binary. These applications are not restricted to discrete data. The early evolutionary trees were all built based on morphological characters such as body size. Moreover, with the burst of new genomic data, such as gene expression, phylogenetic models for continuous traits are again becoming important; see~\cite{hiscott2016efficient} and references therein. Latent tree models in this evolutionary context are also popular in linguistics~\cite{ringe2002indo,shiers2014gaussian}, where a tree represents evolution of languages with modern languages represented by the leaves. The data here typically consists of the acoustic structure of spoken words and are continuous although there are also approaches using syntatic (discrete) data; see, for example,  \cite{guardiano2005parametric,shu2016syntactic}. A related application of latent tree models is in network tomography, where it is used to determine the structure of the connections in the Internet~\cite{castro2004network,eriksson2010toward}. In this application messages are transmitted by sending packets of bits from a source vertex to different destinations and the correlation in arrival times is used in order to infer the underlying network structure. A common assumption is that the underlying network has a tree structure. Then Gaussian latent tree models form a natural correlation model for  arrival times because the correlations diminish with the distance on a tree; see Section~\ref{sec:gaussian}. Typically in this context a special submodel called the \emph{Brownian motion tree model} is used. 

Third, a tree can represent  hierarchical structure in complex data sets. This viewpoint was behind the definition of latent tree models that emerged in the machine learning community~\cite{bishop1998hierarchical,lawrence2004gaussian,zhangCluster}. In a rooted tree every internal vertex represents a cluster given by all the leaves that descent from it. We refer to~\cite{mourad2013survey} for an overview of potential applications. We emphasise that in those applications, it is often not realistic to assume that the true data-generating distribution lies in the latent tree model, which leads to some subtleties in inference; see also Section~\ref{sec:selection}. Another application in this vein is in computer vision and image processing; see~\cite{willsky2002multiresolution} and reference therein. One promising application along these lines is the use of context in computer vision; see, for example,~\cite{choi2010exploiting}.  

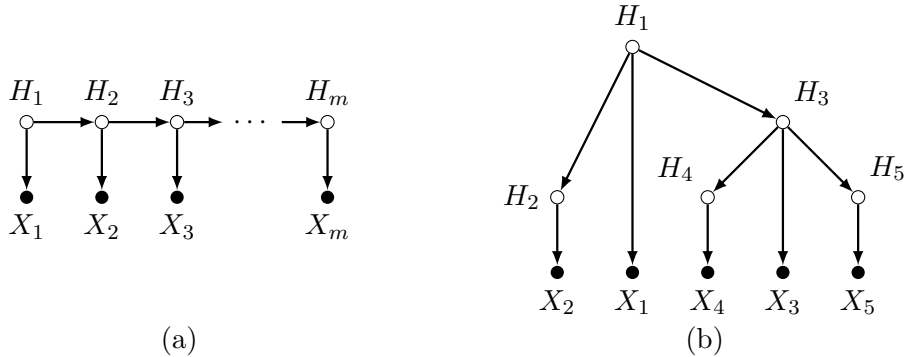
\begin{figure}\centering
	\begin{minipage}{.45\textwidth}\centering
	\tikzstyle{vertex}=[circle,fill=black,minimum size=5pt,inner sep=0pt]
\tikzstyle{hidden}=[circle,draw,minimum size=5pt,inner sep=0pt]
  \begin{tikzpicture}[scale=1]
  \node[vertex] (1) at (0,0)  [label=below:$X_1$] {};
    \node[vertex] (2) at (1,0)  [label=below:$X_2$] {};
    \node[vertex] (3) at (2,0) [label=below:$X_3$]{};
    \node[vertex] (4) at (4,0) [label=below:$X_m$]{};
    \node[hidden] (a) at (0,1)  [label=above:$H_1$]{};
    \node[hidden] (b) at (1,1) [label=above:$H_2$]{};
    \node[hidden] (c) at (2,1) [label=above:$H_3$]{};
     \node (d) at (3,1)  {$\cdots$};
     \node[hidden] (e) at (4,1)  [label=above:$H_m$]{};
           \draw[->,-latex,line width=.3mm] (a) to (b);
           \draw[->,-latex,line width=.3mm] (b) to (c);
           \draw[->,-latex,line width=.3mm] (c) to (d);
           \draw [->,-latex,line width=.3mm](d) to (e);
               \draw[->,-latex,line width=.3mm] (a) to (1);
    \draw[->,-latex,line width=.3mm] (b) to (2);
    \draw[->,-latex,line width=.3mm] (c) to (3);
    \draw[->,-latex,line width=.3mm] (e) to (4);
          \end{tikzpicture}
	\end{minipage}
\begin{minipage}{.45\textwidth}\centering
	\tikzstyle{vertex}=[circle,fill=black,minimum size=5pt,inner sep=0pt]
\tikzstyle{hidden}=[circle,draw,minimum size=5pt,inner sep=0pt]
  \begin{tikzpicture}[scale=1]
  \node[vertex] (1) at (2,1)  [label=below:$X_1$] {};
    \node[vertex] (2) at (1,1)  [label=below:$X_2$] {};
    \node[vertex] (3) at (4,1) [label=below:$X_3$]{};
    \node[vertex] (4) at (3,1) [label=below:$X_4$]{};
    \node[vertex] (5) at (5,1) [label=below:$X_5$]{};
    \node[hidden] (a) at (2,4)  [label=above:$H_1$]{};
    \node[hidden] (b) at (1,2) [label=left:$H_2$]{};
    \node[hidden] (c) at (4,3) [label=85:$H_3$]{};
     \node[hidden] (d) at (3,2)  [label=120:$H_4$]{};
     \node[hidden] (e) at (5,2)  [label=80:$H_5$]{};
           \draw[->,-latex,line width=.3mm] (a) to (b);
           \draw[->,-latex,line width=.3mm] (a) to (c);
           \draw[->,-latex,line width=.3mm] (c) to (d);
           \draw [->,-latex,line width=.3mm](c) to (e);
               \draw[->,-latex,line width=.3mm] (a) to (1);
    \draw[->,-latex,line width=.3mm] (b) to (2);
    \draw[->,-latex,line width=.3mm] (c) to (3);
    \draw[->,-latex,line width=.3mm] (d) to (4);
        \draw[->,-latex,line width=.3mm] (e) to (5);
          \end{tikzpicture}
 \end{minipage}\\
    \begin{minipage}{.45\textwidth}\centering
  (a)
  \end{minipage}  \begin{minipage}{.45\textwidth}\centering
  (b)
  \end{minipage}          \caption{(a) The caterpillar tree defining the hidden Markov model. (b) An example of a tree defining a hidden Markov tree model.}\label{fig:hmm}
\end{figure}

Latent tree models can be also used as a generalization of hidden Markov models. A hidden Markov model (HMM) is a latent tree model on the caterpillar tree in Figure~\ref{fig:hmm}(a). Typically hidden Markov models are \emph{homogeneous} in the sense that the conditional distributions $p(H_i|H_{i-1})$ and $p(X_i|H_i)$ do not depend on $i$; see~\cite{rabiner1989tutorial}. The unobserved part in HMMs follows a Markov chain and so respects a very simple dependency structure. This may be a good approximation of the real dependency structure in the context of time series but is often too restrictive in other applications. Hidden Markov Tree models (HMTMs)~\cite{crouse1998wavelet} relax this restriction allowing for any tree structure on the unobserved vector. Still however, like for HMMs, we assume that each internal vertex has a leaf as a child; see Figure~\ref{fig:hmm}(b). Models of this type were proposed for wavelet-based statistical signal processing but other applications emerged recently: image processing~\cite{choi2001multiscale,romberg2001bayesian}, biomedicine~\cite{makhijani2012accelerated,pfeiffer:CNE23820}, computational lynguistics~\cite{vzabokrtsky2009hidden}. In these settings the unobserved vector typically is binary and the observed part is Gaussian.

This leads to another important reason to study latent tree models: many popular models are submodels of the latent tree models. Examples are given by HMMs, HMTMs, all phylogenetic tree models, but also one factor analysis model, latent class models, and Brownian motion tree models. Models of these types are used virtually everywhere: in biostatistics, machine learning, and social sciences. As noted by Wainwright~and~Jordan in~\cite{wainwright2008graphical}, a more general viewpoint gives a unifying framework for existing algorithms used for these different model classes.

\subsection{Parsimonious latent tree models}\label{sec:statespace}

In this section we briefly discuss certain redundancy in \emph{discrete} latent tree models. We start by discussing the latent class model, which is a special latent tree model where the state space $\cY$ is finite and the underlying tree is a star as in Figure~\ref{fig:6}(a). In this case the tree is typically rooted at the internal vertex and the parameter consists of the root distribution together with the conditional distributions of the leaves given the internal vertex. The number of latent classes corresponds to the number of states $|\cH|$ of the unobserved variable $H$. This relatively simple class of models gives us first insights into general latent tree models. The following result shows that, if the number of unobserved classes is high enough the model becomes \emph{saturated}, that is, it contains all probability distributions over $\cX$. 
\begin{proposition}\label{prop:satur}
If a given latent class model \emph{is not} saturated then
$|\cH|<{|\cX|}/{\max_i |\cX_i|}$.\end{proposition}
The proof of this result can be recovered from the proof of~\cite[Theorem 3]{zhangCluster}.

A latent tree model is \emph{parsimonious} if there is no other latent tree model with a smaller number of parameters that gives the same family of probability distributions over the observed variables. A \emph{discrete} latent class model is \emph{not} parsimonious if $|\cH|\geq  |\cX|/\max_i |\cX_i|$ because, by Proposition~\ref{prop:satur}, every such model is saturated, and the only parsimonious latent tree model over $\cX$ that is saturated is a tree with a single vertex representing the whole vector $(X_1,X_2,X_3)$ as a single variable with $|\cX|$ states.   

A discrete latent tree model is \emph{regular} if the inequality in Proposition \ref{prop:satur} holds for any unobserved vertex $v$ with neighbours $N_v\subset V\setminus \{v\}$, that is, when
\begin{equation}\label{eq:hidclass}
	|\cY_v|<\prod_{u\in N_v}|\cY_u|/\max_u |\cY_u|.
\end{equation} The argument for the latent class model can be generalised to conclude the following; see~\cite{mourad2013survey}: 
\begin{proposition}\label{prop:reducestate}
Any parsimonious discrete latent tree model is regular.\end{proposition}
Proposition~\ref{prop:reducestate} substantially reduces the space of possible discrete latent tree models to consider. As demonstrated by~\cite{zhangCluster} the size of the space of regular models is bounded by $2^{3m^2}$. There are two problems with that in practice. First, this space is still very big with no clear structure. Second, no necessary and sufficient conditions to assure parsimony are known in general and so we do not know how good this reduction is. This leads us to a more tractable subclass of discrete latent tree models called general Markov models, which we discuss in the next section.

\subsection{Gaussian and general Markov models}\label{sec:GMM}

A \emph{general Markov model} is a latent tree model for which all $\cY_v$ are equal and $d=|\cY_v|<\infty$. Models of this type, for $d=2,4,20$ or $61$ appeared in phylogenetics half a century ago~\cite{cavallisforza1967pam}, they were formulated in the most general form over 30 years ago~\cite{barry1987statistical}, but only recently they are becoming increasingly popular; see for example~\cite{allman2008pia,jayaswal2007estimation}. In statistics, general Markov models appeared in the context of causal inference~\cite{pearlBNbook}, or simply as the simplest interesting family of graphical models with unobserved variables~\cite{settimi2000gma}. As we present here, general Markov models stand out as the tractable class of discrete latent tree models, as much as Gaussian models stand out as the tractable class in the continuous case.

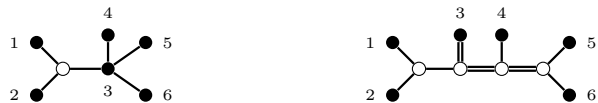
\begin{figure}[htp!]
\centering
		 \tikzstyle{vertex}=[circle,fill=black,minimum size=5pt,inner sep=0pt]
\tikzstyle{hidden}=[circle,draw,minimum size=5pt,inner sep=0pt]
\tikzstyle{every node}=[font=\tiny]
\begin{tikzpicture}[scale=.5]
    \node[vertex] (2) at (.3,-.7) [label=left:$2$]{};
    \node[vertex] (3) at (.3,.7) [label=left:$1$] {};
    \node[vertex] (4) at (2.2,.9) [label=above:$4$]{};
    \node[vertex] (5) at (3.2,0.7) [label=right:$5$]{};
    \node[vertex] (6) at (3.2,-.7) [label=right:$6$]{};
    \node[hidden] (a) at (1,0) {};
    \node[vertex] (c) at (2.2,0) [label=below:$3$] {};
    \draw[line width=.3mm] (a) to (2);
    \draw[line width=.3mm] (a) to (3);
    \draw[line width=.3mm] (a) to (c);
    \draw[line width=.3mm] (c) to (4);
   \draw[line width=.3mm] (c) to (5);
   \draw[line width=.3mm] (c) to (6);
  \end{tikzpicture}\qquad\qquad\qquad
  \begin{tikzpicture}[scale=.5]
  \node[vertex] (1) at (-.7,.7)  [label=left:$1$] {};
    \node[vertex] (2) at (-.7,-.7) [label=left:$2$]{};
    \node[vertex] (3) at (1.1,.9) [label=above:$3$]{};
    \node[vertex] (4) at (2.2,.9) [label=above:$4$]{};
    \node[vertex] (5) at (4,0.7) [label=right:$5$]{};
    \node[vertex] (6) at (4,-.7) [label=right:$6$]{};
    \node[hidden] (a) at (0,0) {};
    \node[hidden] (b) at (1.1,0) {};
    \node[hidden] (c) at (2.2,0) {};
     \node[hidden] (d) at (3.3,0) {};
     \draw[line width=.3mm] (a) to (1);
    \draw[line width=.3mm] (a) to (2);
    \draw[line width=.3mm] (a) to (b);
    \draw[line width=.3mm,double] (b) to (3);
    \draw[line width=.3mm,double] (b) to (c);
    \draw[line width=.3mm] (c) to (4);
   \draw[line width=.3mm,double] (d) to (c);
   \draw[line width=.3mm] (d) to (5);
   \draw[line width=.3mm] (d) to (6);
  \end{tikzpicture}
  \caption{Two equivalent latent tree models. In the model on the right extra restrictions are put on parameters so that the double edges represent equality of random variables.}\label{fig:reduc}
\end{figure}

In the previous section we argued that the space of parsimonious latent tree models is not easy to handle. On the other hand, for general Markov models, inequality~(\ref{eq:hidclass}) is satisfied as long as all internal vertices have degree at least three. For general Markov models and Gaussian latent tree models the necessary and sufficient condition for model parsimony is that:
\begin{itemize}
	\item[(A1)] each unobserved variable has degree at least three,
	\item[(A2)] any two neighbouring variables are neither functionally related nor independent. 
\end{itemize}
It is standard to assume that the underlying tree has no degree two vertices. In fact, for general Markov models and Gaussian latent tree models we can always suppress degree two vertices without changing the model, and so (A1) is always satisfied; see \cite[Section 5.3.4]{LTbook}. In particular, the model defined over a binary rooted tree is equal to the model over the corresponding undirected trivalent tree obtained by suppressing the root. We also always assume the following:
\begin{itemize}
	\item[(A3)] all $Y_v$ are \emph{nondegenerate} meaning that the distribution of $Y_v$ has the full support $\cY_v$.
\end{itemize}
Working without assuming (A2) is sometimes convenient because of the following result, which allows us to focus on learning latent tree models over trivalent trees.
\begin{theorem}\label{th:binaryred}
	Every discrete latent tree model satisfying (A1) and (A2) is a submodel of a latent tree model over a trivalent tree that satisfies only (A1). The same applies to Gaussian latent tree models.
\end{theorem}
A formal proof can be based on the following two observations; see~\cite[Section 5.2.2]{LTbook} for more details. First, a tree with no degree two vertices can be obtained from a trivalent tree by edge contraction. Second, if $Y_u=Y_v$ in a tree model, there exists a simpler model on a tree obtained from $T$ by contracting the edge $(u,v)$. This means that every latent tree model can be realised as a submodel over a trivalent tree with some of the vertices identified. For example, consider the discrete latent tree model for the tree on the left of Figure~\ref{fig:reduc}. Here one of the internal vertices, labeled with $3$, represents an observed random variable. We can alternatively consider a discrete latent tree model on the right of Figure~\ref{fig:reduc}. Here the double edge represent equality between adjacent random variables, and so, the corresponding conditional distribution is degenerate. Directly from the way these models are parameterised, we see that both models are equivalent. 

There are numerous advantages of general Markov models. The link to statistical physics, phylogenetics, and Markov processes allows to develop efficient algorithms with strong theoretical guarantees. General Markov models retain very good performance on real-world data~\cite{choi:2011}. The space of parsimonious general Markov models is also much more tractable than the space of all regular latent tree models over $\cX$. In some applications we cannot assume that the state spaces of observed variables are equal but we can assume that the unobserved variables have the same state spaces. This holds for HMMs and HMTMs. Many techniques discussed later can be generalised to that case.

For any edge $u\to v$ denote by $M^{uv}$ the $d\times d$ stochastic matrix representing the conditional distribution $p_{v|u}$. In analogy to Markov chains, $M^{uv}$ is called a \emph{transition matrix} and each of its \emph{rows} represents a conditional distribution of $Y_v$ given a particular value of $Y_u$. In the general Markov model, each $M^{uv}$ is an arbitrary stochastic matrix. Models used in phylogenetics are usually more constrained. They are generated through a continuous time  Markov process on $T$ with a given \emph{rate matrix} $Q$, that is a real $d\times d$ matrix with row elements summing to zero and all off-diagonal elements nonnegative. In this setting the transition matrices are given by   
\begin{equation}\label{eq:Mexp}
M^{uv}\;=\;\exp(t_{uv}Q),	
\end{equation}
where $t_{uv}>0$ is an edge parameter and $\exp(\cdot)$ is the matrix exponential function. The rate matrix $Q$ typically has some further symmetries. For example, in the Jukes-Cantor model~\cite{jukes1969epm}, all off-diagonal entries of $Q$ are equal; see~\cite{felsenstein2004inferring} for a review of most popular phylogenetic models. 

%Phylogenetic models are submodels of general Markov models. To understand the relation between these two models we need to understand better the additional positivity conditions imposed by $M^{uv}$ by the exponential representation in~(\ref{eq:Mexp}). A classical result by Kingman~\cite{kingman1962imbedding} shows that, if $d=2$, $M^{uv}=\exp(\tau_{uv}Q)$ for some rate matrix $Q$ and $\tau_{uv}>0$ if and only if $\det M^{uv}>0$. If $d>2$, the determinant of $M^{uv}$ must be positive by the Jacobi identity but no necessary and sufficient conditions are known. 
%
%

\section{Second-order moment structure}

Many learning algorithms for latent tree models use the second-order structure of the observed distribution, that is, correlations between the observed variables, mutual information, or other aspects of the pairwise marginal distributions. These algorithms typically exploit links to tree metrics. In this section we start by explaining this link and how it can be used to design robust learning algorithms.

\subsection{Gaussian latent tree model}\label{sec:gaussian}

The earliest example of Gaussian latent tree models is the factor analysis model with a single unobserved factor~\cite{anderson1956statistical,thurstone1934vectors}. More general Gaussian tree models were not studied until more recently~\cite{choi:2011,ASSZ2014}. In the fully observed Gaussian tree model the \emph{inverse covariance} matrix of $Y$ is very sparse. However, the inverse covariance matrix of the observed subvector $X$ has no zeros and it is not convenient to work with. On the other hand, the \emph{covariance} matrix of $X$ provides a great insight into the structure of Gaussian latent tree models and latent tree models in general. 

Consider any two vertices $i,j$ of $T$. Using the Markov properties implied by the tree model, the correlation $\rho_{ij}={\rm corr}(Y_i,Y_j)$ can be written as the product (see, e.g.,~\cite{choi:2011})
\begin{equation}\label{eq:RhoPath}
	\rho_{ij}\;=\;\prod_{(u,v)\in\overline{ij}} \rho_{uv},
\end{equation} 
where $\overline{ij}$ denotes the unique path between $i$ and $j$ in $T$. Restricting~(\ref{eq:RhoPath}) only to pairs of leaves $i$, $j$ gives the parameterization of the correlations of the Gaussian latent tree model. In particular, vertices that are far from each other in the tree tend to be less correlated. This appealing property makes this model useful for hierarchical data clustering and network tomography problems as described in Section~\ref{sec:apps}.

The constraints on the correlations induced by~(\ref{eq:RhoPath}) are the only nontrivial constraints on Gaussian latent tree models. The variances of the observed variables can be arbitrary and together with the \emph{edge correlations}, $\rho_{uv}$ for all edges $(u,v)$, they provide a set of parameters for the latent tree model. The variances of the unobserved variables do not affect the observed distribution and so typically are set to $1$. The model is parsimonious as long as it satisfies (A1) and (A2). The condition (A2) is satisfied if the edge correlations satisfy $|\rho_{uv}|\in (0,1)$.

\begin{example}\label{ex:Gquartet}
	Consider the quartet tree in Figure~\ref{fig:6}(b) and denote the internal vertices by $u$, $v$, where $u$ is closer to $1$ and $2$. If correlations between the four leaves come from a latent tree model on this quartet tree, then there is a collection of edge correlations $\rho_{1u}$, $\rho_{2u}$, $\rho_{uv}$, $\rho_{3v}$, and $\rho_{4v}$ all with values in $[-1,1]$ such that $\rho_{12}=\rho_{1u}\rho_{2u}$, $\rho_{13}=\rho_{1u}\rho_{uv}\rho_{3v}$, $\rho_{14}=\rho_{1u}\rho_{uv}\rho_{4v}$, $\rho_{23}=\rho_{2u}\rho_{uv}\rho_{3v}$, $\rho_{24}=\rho_{2u}\rho_{uv}\rho_{4v}$, and $\rho_{34}=\rho_{3v}\rho_{4v}$. In particular
	\begin{equation}\label{eq:quarteteqs}
	\rho_{13}\rho_{24}\;=\;\rho_{14}\rho_{23}\;=\;\rho_{12}\rho_{34}\rho_{uv}^2.	
	\end{equation}
	The first equality is an example of a nontrivial relation between the observed correlations. The second equality can be used to recover the value of the edge parameter $\rho_{uv}$ given only the observed correlations, or to establish the inequality $\rho_{14}\rho_{23}\geq \rho_{12}\rho_{34}$. \end{example}
A systematic study of polynomial constraints defining statistical models is part of algebraic statistics \cite{oberwolfach2009}. For Gaussian latent tree models these constraints were studied in \cite[Section 6]{sullivant2008agg} and \cite{ASSZ2014}.

Correlation matrices of latent Gaussian tree model are closely linked to tree metrics. Given a tree $T$ with $m$ leaves assign to each of each edges $u-v$ a nonnegative length $d_{uv}$. With this choice we can now compute the distance between any two leaves $i,j$ summing the lengths of the edges on the unique path between $i$ and $j$, that is,  
\begin{equation}\label{eq:treemet}
	d_{ij}\;\;=\;\;\sum_{(u,v)\in \overline{ij}} d_{uv}.
\end{equation}  
Consider an $m\times m$ symmetric matrix $D=[d_{ij}]$ with zeros on the diagonal. Call $D$ a \emph{tree metric} if its entries satisfy~(\ref{eq:treemet}) for some tree $T$ and edge lengths.

%\begin{definition}
%The space of \emph{tree metrics} is the set of $m\times m$ symmetric matrices $D=[d_{ij}]$ such that, for some tree $T$ with $m$ leaves and a collection of edge lengths $d_{uv}>0$ for edges $(u,v)$ of $T$, the entries of $D$ satisfy $d_{ij}=\sum_{(u,v)\in \overline{ij}} d_{uv}$.
%\end{definition}

Let $\Sigma=[\rho_{ij}]$ be a correlation matrix in a Gaussian latent tree model. Assume first that $\rho_{ij}\neq 0$ for all $1\leq i<j\leq m$. Consider a symmetric matrix $D=[d_{ij}]$ with $d_{ij}=-\log |\rho_{ij}|$. Then~(\ref{eq:RhoPath}) translates into~(\ref{eq:treemet}). Since $|\rho_{uv}|\in (0,1]$ for all $(u,v)\in \overline{ij}$ also $d_{uv}\geq 0$, and so $D$ is a tree metric. If $\Sigma$ contains zero entries~(\ref{eq:RhoPath}) implies that zeros cannot appear arbitrarily. For every three indices $i,j,k$ if $\rho_{ij}\neq 0$ and $\rho_{jk}\neq 0$ then also $\rho_{ik}\neq 0$. It follows that the correlation matrix in a Gaussian latent tree model has a block diagonal structure and within each block all entries are non-zero. Each block can be transformed to a tree metric on the corresponding subtree. 

The connection to tree metrics will be exploited in many ways. For example, it is a standard result in phylogenetics~\cite{buneman1971} that for any tree metric $D$ the underlying tree and positive edge lengths $d_{uv}$ can be recovered uniquely. Because $|\rho_{uv}|=\exp(-d_{uv})$, we obtain the following corollary.
\begin{theorem}\label{th:identifyGaussian}
If $\Sigma=[\rho_{ij}]$ is a correlation matrix in a Gaussian latent tree model with non-zero entries, then the underlying tree and edge correlations are identified uniquely (up to sign). 
\end{theorem}
To have a concrete example of how it works consider again the quartet tree model in Example~\ref{ex:Gquartet}. Suppose that we are given a correlation matrix from the Gaussian latent tree model but we do not know the underlying tree. By Theorem~\ref{th:binaryred} we can first constrain ourselves to trivalent trees, that is, the three possible quartets: $12/34$, $13/24$, and $14/23$, where this notation indicates how leaves group together. In the first case, by~(\ref{eq:quarteteqs}), we have 
\begin{equation}\label{eq:quartetseqs}
|\rho_{12}\rho_{34}|\;\geq\; |\rho_{13}\rho_{24}|\;=\;|\rho_{14}\rho_{23}|	
\end{equation}
and by symmetry we obtain similar relations for the other two trees. Therefore, we can find the underlying tree by computing quantities of the form $|\rho_{ij}\rho_{kl}|$ and choosing the largest. If they all happen to be equal, the underlying tree is the star tree. 

We can extend Theorem~\ref{th:identifyGaussian} to arbitrary correlation matrices but here identifiability of the edge correlations is more involved; see~\cite{DLWZ2014} for details. For instance, consider again the model in Example~\ref{ex:Gquartet}. Suppose that $\rho_{ij}=0$ for all $i\in \{1,2\}$ and $j\in \{3,4\}$ and $\rho_{12},\rho_{34}\neq 0$. Then we can easily identify the underlying tree $12/34$. Indeed, for no other tree with four leaves there is a choice of edge correlations giving precisely this pattern of zeros. Identifying the parameters is harder. We check that $\rho_{uv}$ must be zero. The other edge correlations must be nonzero but they are identified only up to the relations $\rho_{12}=\rho_{1u}\rho_{2u}$ and $\rho_{34}=\rho_{3v}\rho_{4v}$.

%In the field of phylogenetics many methods have been developed to learn the underlying tree and edge lengths from observed leaf distances. These methods can be used also to learn the underlying tree in the Gaussian latent tree model together with the underlying edge correlations $\rho_{uv}$. These methods are all provably consistent if the original tree-learning methods are. We discuss this in more detail in Section~\ref{sec:distance}.

\subsection{General Markov models}\label{sec:discrete}

Tree metrics appear also in the description of discrete latent tree models \cite{barry1987statistical,lake1994ret,lockhart1994recovering,steel2004}. For any edge $u\to v$ denote by $M^{uv}$ the $d\times d$ transition matrix representing the conditional distribution $p_{v|u}$. Denote by $P^{uv}$ the $d\times d$ matrix of the marginal distribution of $(Y_u,Y_v)$, and by $P^{uu}$ a diagonal matrix with the marginal distribution of $Y_u$ on the diagonal. For any two vertices $u,v$ let 
\begin{equation}\label{eq:u}
\tau_{uv}\;:=\;\frac{\det(P^{uv})}{\sqrt{\det(P^{uu}P^{vv})}},\end{equation}
where the denominator is non-zero by assumption (A3). By essentially the same argument as in~\cite[Theorem 8.4.3]{semple2003pol} we obtain the following path-product formula
\begin{equation}\label{eq:qij}
\tau_{ij}\;\;\;=\;\;\prod_{(u,v)\in \overline{ij}} \tau_{uv}\qquad\mbox{for all }i,j\in V.
\end{equation}
%where 
%$$
%\tau_{uv}=\det M^{uv}\sqrt{\frac{\det (P^{uu})  }{\det(P^{vv})}}.
%$$
In the case of binary variables, $\det P^{ij}={\rm cov}(X_i,X_j)$, $\det(P^{ii})={\rm var}(X_i)$ and so $\tau_{ij}$ is the correlation, which implies that~(\ref{eq:qij}) reduces to~(\ref{eq:RhoPath}). 
%\begin{remark}
%We obtain the path-product formula for correlations like in~(\ref{eq:RhoPath}) as long as $H$ is binary or Gaussian irrespective of the form of $X$. In particular it holds for the hidden Markov tree model with binary unobserved variables. 
%\end{remark}

In general the interpretation of the edge parameters $\tau_{uv}$ is more complicated. Using the identity $P^{uu}M^{uv}=P^{uv}$ we can write 
\begin{equation}\label{eq:u2}
\tau_{uv}\;\;=\;\;\det M^{uv}\sqrt{\frac{\det (P^{uu})  }{\det(P^{vv})}}.
\end{equation} 
Therefore, we have $\tau_{uv}=0$ if and only if $\det M^{uv}=0$. If $d=2$ this is equivalent to independence of $Y_u$ and $Y_v$ but in general it is a strictly weaker condition. For example, if $d=3$, and the first two rows of $M^{uv}$ are equal but not equal to the third one, then $\det M^{uv}=0$ but $Y_u$ and $Y_v$ are not independent. Exactly like in the Gaussian case the edge parameters satisfy  $|\tau_{uv}|\leq 1$ and the border values $\pm 1$ consider to functional dependence of $Y_u$ and $Y_v$. Indeed, by the Bayes' theorem, we have
	$$
	M^{vu}\;=\; (P^{vv})^{-1}(M^{uv})^TP^{uu}.
	$$
	With this notation,~(\ref{eq:u2}) gives
	$$
	\tau_{uv}^2=\det (M^{uv}M^{uv}) \det(P^{uu}(P^{vv})^{-1})=\det M^{uv} \det M^{vu}. 
	$$
	Because both $M^{uv}$ and $M^{vu}$ are stochastic matrices, all their eigenvalues lie in the unit circle. In particular, $\tau_{uv}\in [-1,1]$ and it is equal to $\pm 1$ precisely when $M^{uv}$ is a permutation matrix, or in other words, if $Y_u$ and $Y_v$ are functionally related.

This again gives a direct link to tree metrics and the following result; see \cite{chang1991reconstruction,steel1993invertible}.
\begin{theorem}
	For any distribution in a general Markov model the underlying tree can be uniquely recovered given only its 2-way margins.
\end{theorem}
We implicitly assume that $\tau_{ij}$ are nonzero but the theorem can be extended to forests. 
%Unlike in the Gaussian case, we cannot typically identify the parameters of the model using only 2-way margins. In Theorem~\ref{th:identify} we present that 3-way margins are always sufficient.
%

\subsection{Linear models} The fact that in the Gaussian and in the binary case the correlation between observed variables decompose as in~(\ref{eq:RhoPath}) can be proved using the following fact: if $X,Y$ are binary (or jointly Gaussian) then the conditional expectation $\E[X|Y]$ is a linear function of $Y$. As we see in this section the discrete case can be also interpreted in this way, which gives a unifying framework to understand~(\ref{eq:RhoPath}) and~(\ref{eq:qij}). Moreover, it gives a much larger families of potential models to consider that also admit a path-product formula linking it to tree metrics. 

Let each variable $Y_u$ in the system be modeled as a random vector in $\R^k$ for a fixed $k$. A ternary variable, for example, will take values $(0,0)$, $(1,0)$, $(0,1)$ in $\R^2$ instead of the typical $1,2,3$ in $\R$.  Each variable can be either discrete or continuous but we add a minor requirement that the matrix $\Sigma_{vv}=\E Y_v Y_v^T-\E Y_v (\E Y_v)^T$ is positive definite. A \emph{linear latent tree model} is a latent tree model in which for every edge $u-w$ in the tree, the conditional expectation $\E[Y_u|Y_w]$ is an affine function of $Y_w$. Linear latent tree models include (multivariate) Gaussian latent tree models, Kalman filters, Gaussian mixtures, Poisson mixtures and general Markov models. Models of this type were first discussed in~\cite{anandkumar2011spectral}. Here we propose a slightly different exposition using ideas from~\cite[Section 4.3]{LTbook}.

We define the normalized version $\bar Y_v$ of $Y_v$ as $\bar Y_v:=(\Sigma_{vv})^{-1/2}(Y_v-\E Y_v)$.  Denoting $\Sigma_{uv}=\E Y_u Y_v^T-\E Y_u (\E Y_v)^T$ we obtain
\begin{equation}\label{eq:condexp}
	\E[\bar Y_u|Y_w]\;\;=\;\;\Sigma_{uu}^{-1/2}\Sigma_{uw}\Sigma_{ww}^{-1/2}\,\bar Y_w.
\end{equation}
 Define $\tau_{uv}:=\det(\Sigma_{uu}^{-1/2}\Sigma_{uv}\Sigma_{vv}^{-1/2})$. By the law of total expectation, it follows from (\ref{eq:condexp}) that $\tau_{uv}=\det (\E[\bar Y_u \bar Y_w^T])$. Let $Y_u,Y_v,Y_w$ be three random variables with values in $\R^k$ such that $Y_u\indep Y_w|Y_u$. Then 
$$
\E[\bar Y_u \bar Y_w^T]=\E\left[\E[\bar Y_u|Y_w](\E[ \bar Y_w^T|Y_w])^T\right]=\Sigma_{uu}^{-1/2}\Sigma_{uw}\Sigma_{ww}^{-1}\Sigma_{wv}\Sigma_{vv}^{-1/2},
$$
which implies that $\tau_{uv}=\det (\E[\bar Y_u \bar Y_w^T])=\tau_{uw}\tau_{wv}$. Applying this argument recursively we  conclude that the path-product decomposition of $\tau_{ij}$ given in~(\ref{eq:qij}) holds for any linear latent tree model. 

This clearly generalizes the Gaussian case. To see that this also generalises~(\ref{eq:qij}), for each each discrete random variable with $d$ states take $k=d-1$ and set the state-space to be $\{0,e_1,\ldots,e_{d-1}\}$, where $0$ is the origin and $e_i$ are the elements of the standard basis of $\R^{d-1}$. 
\begin{proposition}
With the above convention $\det(\Sigma_{uu})=\det(P^{uu})$ and $\det(\Sigma_{uv})=\det{P^{uv}}$.	\end{proposition}
The proposition implies that $\tau_{uv}$ as defined in this section is equal to $\tau_{uv}$ as defined in Section \ref{sec:discrete}.
\begin{proof}[Proof of the proposition]
	Consider the matrix $A$ obtained from $P^{uv}$ by elementary row and column operations: add all rows to the first row and all columns to the first column. Basic linear algebra implies $\det P^{uv}=\det A$. Matrix $A$ has the following block structure. The top-left  $1\times 1$-block is equal to $1$. The bottom-right $(d-1)\times (d-1)$-block is equal to $\E Y_u Y_v^T$ and the remaining two blocks are $\E Y_u$ and $\E Y_v^T$. The formula for the determinant of a block matrix implies that $\det A=\det(\E Y_u Y_v^T-\E Y_u\E Y_v^T)=\det \Sigma_{uv}$. The proof of the other equality is analogous.  
\end{proof}

%Learning based on this representation can be done in a similar way as discussed, for the Gaussian case, in the end of Section~\ref{sec:gaussian}. A straightforward approach is to try to learn the underlying quartets by using constraints of the form~(\ref{eq:quartetseqs}). In the later sections we will briefly discuss some guarantees for such a method. For more details see~\cite{anandkumar2011spectral}. In the next section we discuss a more global approach.

\subsection{Distance based methods}\label{sec:distance}

The maximum likelihood tree topology recovery is NP hard~\cite{roch2006short}. This has motivated a number of investigations of other tractable methods for learning trees as well as theoretical guarantees on performance. The link between tree metrics and latent tree models described in the previous sections makes it possible to come up with consistent methods to learn a tree that work in polynomial time. This approach dates back to~\cite{barry1987statistical}; see~\cite{holland01012013} and references therein. 

For a concrete example consider a sample from the Gaussian latent tree model. Given the sample correlation matrix with elements $\hat\rho_{ij}$ we compute distances $\hat d_{ij}=-\log|\hat\rho_{ij}|$. Now use any of the methods to learn a tree metric from observed distances. This gives a tree $\hat T$ and edge distances $\hat d_{uv}$, or equivalently absolute values of the edge correlations $\hat\rho_{uv}$. Such a method will be (statistically) consistent given the original tree distance method is \emph{consistent}, which in this context means that the method outputs the correct tree given a tree metric.

There are many methods that try to recover the underlying tree from noisy distances. The most popular are the Neighbour-Joining (NJ) algorithm and the least-squares method but many other algorithms are available; see Section 7.3 in~\cite{semple2003pol} for an overview. All popular methods are both well studied and widely implemented, for example, in \textsc{R}; see Section 5.1 in~\cite{paradis2011analysis}. Most of the methods, including NJ and the least squares, are {consistent}. We also note in passing that these methods output an undirected tree but rooting is also possible by finding an appropriate outgroup; see ~\cite[Section 7.3]{durbin1998bsa}.

An appealing property of this method, as applied in the Gaussian case, is that there is a one-to-one correspondence between edge lengths and model parameters given by edge correlations (up to sign). In the discrete case the situation is more complicated. Here we obtain noisy distances by defining $\hat \tau_{ij}$ like in~(\ref{eq:u}) with $P^{ij}$, $P^{ii}$, and $P^{jj}$ replaced by their sample versions. This gives $\hat d_{ij}=-\log|\hat \tau_{ij}|$. Using the NJ algorithm we obtain an estimate of the underlying tree and parameters $\tau_{uv}$. However, this is in general not enough to recover all model parameters. A special case when it is possible is so called \emph{symmetric discrete distributions}. In these submodels the matrix $M^{uv}$ has all off-diagonal entries equal and the root distribution is assumed to be uniform. In statistical mechanics this corresponds to the Potts model, which in the binary case gives the Ising model; see Example 3.2 in \cite{wainwright2008graphical}.

Although using the NJ algorithm as a tree learning subroutine results in a computationally efficient and consistent method, consistency is only a minor requirement, and other tree learning procedures can be preferred in order to allow for a more sophisticated statistical analysis. An early example is the Dyadic Closure Tree Construction method (DCTC)  \cite{erdos1999few} and witness-antiwitness method (WAM) \cite{erdos199977} both focusing on learning the underlying tree by learning certain quartets that hold; see the latter paper for more details and references. In~\cite{choi:2011} the authors propose two other algorithms: recursive grouping (RG) and CLGrouping. Recursive grouping builds the latent tree recursively by identifying sibling groups using distances $d_{ij}$. CLGrouping starts with a pre-processing procedure in which a tree over the observed variables is constructed, or more precisely, the Chow--Liu tree. This global step groups the observed vertices that are likely to be close to each other in the true latent tree, thereby guiding subsequent recursive grouping (or equivalent procedures such as neighbour-joining) on much smaller subsets of variables. This results in more accurate and efficient learning of latent trees. This can be further improved by using distance information to learn locally small trees and then glue them together. An example of such a divide-and-conquer algorithm is given in~\cite{huang2014scalable}.

The distance based methods to learn the underlying tree are based predominantly on second order margins and so typically are very robust with respect to the sampling error and model misspecifications. In the next section we present other methods to estimate model parameters that use much more information about the underlying distribution.

%An advantage of these distance based methods is that they can be easily generalised to nonparametric contexts. For example, ~\cite{song2011kernel} proposed kernel based methods that perform well under a very general scenario.\piotr{FINISH}
%

\section{Selected theoretical results}\label{sec:theory}

Latent tree models, like all models with unobserved variables, suffer from various problems and learning is generally complicated. The complex geometry of models with unobserved variables usually leads to difficulties in establishing the identifiability of their parameters, and the likelihood function has many local maxima, which lie on the boundary of the parameter space; see, for example,~\cite{chor2000mml,pwz-2009-semialgebraictrees}. In consequence, standard inference and model selection procedures are not fully justified in this setting. In this section we discuss parameter identifiability, sample complexity, and model selection methods. These three seemingly unrelated topics all deliver one important message: the class of latent tree models is well behaved from the statistical point of view as long as all the edge correlations are sufficiently large in the absolute value. Otherwise, it should be used with caution. Some further issues with the likelihood function for this model class will be discussed in Section \ref{sec:like}.

\subsection{Identifiability}\label{sec:ident}

A parametric model $(P_\theta)$ is \emph{identifiable} if $P_\theta=P_{\theta'}$ implies $\theta=\theta'$. In other words, the parameterization map $\theta\mapsto P_\theta$ is a bijection between the parameter space and the model. Latent tree models are never identifiable. This can be seen for latent class models, where permuting labels of the unobserved variable makes no difference in the observed distribution; see, for example,~\cite{zou06042011}. This is known as the \emph{label swapping problem}. The label swapping is not a serious problem in practice. We can always take account of it by restricting the parameter space. However, this still does not make the model identifiable because there are special subspaces in the parameter space that map to the same observed distribution. We illustrate this issue in the simplest possible example.

\begin{example}\label{ex:3gaussian}
Consider the Gaussian latent tree model on the star tree with three leaves, see Figure~\ref{fig:6}(a). Denote $\rho_{ij}={\rm corr}(X_i,X_j)$ and $\rho_i={\rm corr}(X_i,H)$ for $i=1,2,3$. By~(\ref{eq:RhoPath}), this model is parameterized by
$\rho_{12}=\rho_1\rho_2$, $\rho_{13}=\rho_1\rho_3$, and $\rho_{23}=\rho_2\rho_3$.  If the observed correlations are non-zero we can identify $\rho_1$ up to the sign and then the other parameters as follows
$$
\rho_1^2=\frac{\rho_{12}\rho_{13}}{\rho_{23}},\qquad \rho_2=\frac{\rho_{12}}{\rho_1},\qquad \rho_3=\frac{\rho_{13}}{\rho_1}. 
$$
Suppose now that some observed correlations vanish. The form of the parameterization implies that it is impossible for only one of them to vanish. If two correlations are zero, say $\rho_{12}=\rho_{13}=0$, then the set of all triples $(\rho_1,\rho_2,\rho_3)$ mapping to $(0,0,\rho_{23})$ is a smooth one-dimensional subset given by $\rho_1=0$ and $\rho_2\rho_3=\rho_{23}$. Suppose now that \emph{all} observed correlations are zero. Then the corresponding parameters form  the union of three intervals
	$$
	\{\rho_1=\rho_2=0\}\cup \{\rho_1=\rho_3=0\}\cup \{\rho_2=\rho_3=0\}.
	$$ \end{example}
This example motivates the following definition.

\begin{definition}
	A model is \emph{generically identifiable}, if the parameterization map is finite-to-one everywhere outside of a measure zero set.
\end{definition} 
Geometrically, generic identifiability means that for a typical point in the model, its preimage under the parameterization map, also called a \emph{fiber}, is a finite collection of points. Showing that the Gaussian latent tree model satisfying (A1) is generically identifiable can be done by arguments as in Example \ref{ex:3gaussian}. For general Markov models generic identifiability is more subtle and, in general, the second-order moments contain not enough information about the underlying parameters. It turns out that three-way margins are already enough. The following is the main result of~\cite{chang1996frm}.
\begin{theorem}\label{th:identify}
	Every general Markov model satisfying (A1) is generically identifiable. In fact, to identify the parameters it is enough to know the $3$-way marginal distributions. 
\end{theorem}
Theorem~\ref{th:identify} does not cover the case when the state-spaces $\cY_v$ are allowed to vary. In this case techniques of~\cite{allman2009identifiability} may be useful to establish identifiability.	

The analysis of when exactly identifiability fails can be in general complicated. In the Gaussian and the binary case these special points correspond to some correlations being zero~\cite{pwz-2010-identifiability}. As illustrated by Example~\ref{ex:3gaussian}, depending on the situation the corresponding fiber can be either a smooth subset of the parameter space or it can be singular. These theoretical results and examples provide the following insight. If the true data-generating distribution is characterised by high correlations between variables, it is also far from any of the special singular points. However, if some variables have a low degree of dependence, then estimation may become difficult and standard asymptotic theory breaks down. 

\subsection{Guarantees for tree reconstruction}\label{sec:sample}

 A basic question regarding learning of latent tree models is that of the \emph{sample complexity} of the tree reconstruction problem. Given an estimator $\hat  T$ of the true tree $T$ we want to assure that $\P(\hat T=T)>1-\delta$ for some fixed small $\delta$. This is only possible if the sample size is big enough. It is known that irrespective of the method $n=(\log m)$ is necessary~\cite{erdos1999few,mossel2003impossibility} but it is typically not enough. The first systematic study of when the logarithmic bound is sufficient was offered in \cite{erdos1999few} for the binary symmetric model and in \cite{erdos199977} for general Markov models. The link to tree metrics shows that, in order to recover the underlying tree with high probability, the edge lengths cannot be too small or too large. For linear latent tree models this means that the edge parameters $\tau_{uv}$ must satisfy $c\leq |\tau_{uv}|\leq C$ for some $0<c<C<1$. Sample bounds we discuss will always necessarily depend on the number of observed variables $m$ and the parameters $c,C$.

% Moreover, in order for the problem to be well-posed we need to restrict to recovering trivalent trees. Otherwise, by Theorem \ref{th:binaryred}, every non-trivalent tree model could be arbitrary well approximated by a trivalent tree model. 

One of the important concepts developed in \cite{erdos1999few,erdos199977} is that of the tree depth: the \emph{depth} of a tree $T$ with vertices $V$ and leaves $W\subset V$ is $\max_{i\in W\setminus V}\min_{j\in V}|\overline{ij}|.$ For example, the depth of the tree in Figure \ref{fig:6}(c) is $2$ and of both trees in Figure~\ref{fig:hmm} is $1$. Latent tree models for trees with few unobserved variables or small depth are generally easier to learn. This leads to one group of results assuming that the depth of the underlying tree is $O(1)$, that is, constant in the number of leaves of $T$. It was shown in \cite[Theorem 14]{erdos199977} that, under constant tree depth, there is an algorithm for general Markov models whose sample complexity is logarithmic in $m$. This has been generalised to the Gaussian case in~\cite{choi:2011}. 

%Similarly strong results are possible without assuming constant depth. The potential problem in that case is that low correlations (large distances) require relatively many samples to estimate. As noted in~\cite{bhamidi2010network}, in this case we can focus on large correlations, build smaller subtrees, and try to recover the underlying tree using only those. This can be made into a procedure that recovers the correct tree with high probability with sample complexity logarithmic in $m$. In phylogenetics an example of such an approach is given in~\cite{daskalakis2011phylogenies,erdos1999few,erdos199977}. This was further successfully applied for general linear latent tree models~\cite{huang2014scalable}.

Similarly strong results are possible without assuming constant depth. It was first conjectured by Mike Steel \cite{steel2001my} that for the binary symmetric latent tree model the sample complexity is logarithmic in $m$ as long as $c\leq |\tau_{uv}|\leq C$
for some $\frac{\sqrt{2}}{2}<c<C<1$, that is, when the parameters lie in the \emph{Kesten--Stigum (KS) regime}. This conjecture was proved in \cite{daskalakis2011phylogenies,mossel2004phase}. In the subsequent work, it has been shown in multiple scenarios, including the Gaussian case, that in the KS regime, high probability tree topology reconstruction may be achieved with $n=O(\log m)$ samples; see~\cite{mossel:2013} and references therein. The KS regime plays also an important role in the sample complexity analysis of the maximum likelihood estimator. In general the sample complexity is polynomial in $m$ under very general assumptions. However, for symmetric models and under a discretisation assumption, it becomes logarithmic in $m$ if the true distribution lies in the KS regime~\cite{roch2015phase}.

\subsection{Model selection}\label{sec:selection}

As presented in the previous section, a lot of effort in the literature is put to obtain performance guarantees for the proposed learning algorithms. This analysis is done under assumption that the data come from a latent tree model. Although this assumption seems to be reasonable for phylogenetics and several other applications mentioned here, it is certainly not so for many other kinds of data (e.g. survey data from medicine and social sciences) for which these models are used, c.f.~Section \ref{sec:apps}. In all these cases model selection techniques tend to outperform algorithms that aim at finding the ``true'' tree~\cite{zhang2017latent}. 

An example of a model selection algorithm is EAST (Expansion, Adjustment and Simplification until Termination)~\cite{chen2012model}, which aims to find the model with the highest Bayesian Information Criterion (BIC) score. Given a random sample $x^{(1)},\ldots,x^{(n)}$ the BIC is \\[-.2cm]
$$
{\rm BIC}(M(T,\cY))\;=\;\ell(\hat\theta,T)-\frac{d}{2}\log n,\\[-.1cm]
$$ 
where $\hat\theta$ is the maximum likelihood estimator, and $d$ is the dimension of the model $M(T,\cY)$. Computing the dimension of a model with unobserved variables is not always easy because it need not correspond to the number of parameters in the model. In that case a more detailed study of the generic rank of the Jacobian of the parameterization is needed. However, for general Markov models, Theorem~\ref{th:identify} can be used to show that simple parameter count does suffice to compute the dimension of the model whenever (A1) holds. For Gaussian latent tree models the analysis in Section~\ref{sec:gaussian} showed that the dimension is $m+|E|$, where $E$ is the set of edges of the underlying tree and $m$ is the number of its leaves.  

Theoretical importance of the BIC criterion~\cite{schwarz1978edm} is that it provides an asymptotic approximation to the marginal likelihood of the model, and so it can guide model selection in the Bayesian setting.  However, latent tree models lead to a new difficulty in that the Fisher information matrix of a latent tree model is singular along certain submodels. Such singularities invalidate the mathematical arguments that lead to the Bayesian information criterion; see, for example,~\cite{drton2013sBIC,shaowei_rlct,watanabe_book}. This again shows that  BIC must be used with caution in the case of week correlations. In that case the BIC may be too conservative and select too small models. Correcting BIC to account for singularities involves an indepth study of the model  geometry. For binary latent tree models this has been done in~\cite{pwz-2010-bic} and for Gaussian latent tree models in~\cite{DLWZ2014}. In general, a simpler adjustment can be done by replacing BIC with Watanabe's WAIC \cite{watanabe2013widely}.

\section{Estimation and inference}\label{sec:learning}

In Section~\ref{sec:distance} we described some natural tree-metric based approaches to learn latent tree models. In this section we briefly describe other popular learning algorithms.

% In the previous section we described a family of learning methods based on the tree metric connection. In this section we briefly describe other popular algorithms. A more theoretical analysis is provided in Section~\ref{sec:theory}.

\subsection{Fixed tree structure}\label{sec:like}

%The distance based methods to learn the underlying tree discussed in Section~\ref{sec:distance} can be used to learn the continuous parameter only in limited number of cases: the Gaussian and the symmetric discrete. In all other cases, at least in principle, the learning task (L2) need to be treated separately. In this section we briefly describe generic estimation methods. In the next section we discuss the structural EM algorithm, which is more specific to latent tree models. 

The maximum likelihood estimator is one of the most popular estimators of the continuous parameter. It cannot be computed in a closed form and the likelihood function is very complicated to analyse directly; see, for example, \cite{chor2000mml,pwz-2009-semialgebraictrees} . However, there are various numerical methods to maximise the likelihood function that were developed in the context of latent tree models or, more generally, Bayesian networks with unobserved variables. The most natural choice is the EM algorithm, which  is easy to set up; see, for example, Section 3.1.1 in~\cite{mourad2013survey}. There are also several methods building upon the idea of the EM algorithm. For example the progressive EM algorithm of~\cite{chen2015progressive} uses a method of moments estimator as a subroutine, which leads to a more computationally efficient estimator. 

As with all other EM-based methods, these approaches depend on the initialisation and suffer from the possibility of being trapped in local optima. This algorithm seems to work well in the case when the observed variables are highly correlated~\cite{wang2006severity}. The situation becomes much more complicated  in the presence of weak correlations as the numerical procedures become unstable. Another problem is that the maximum likelihood estimator  often lies on the boundary of the parameter space~\cite{pwz-2009-semialgebraictrees}, which has many important consequences. First, the gradient of the likelihood function does not vanish at such a point and so the standard asymptotics does not apply. Second, there may be relatively distant points in the parameter space that give a similar value of the likelihood function as the maximum likelihood estimator. Finally, the boundary points typically correspond to situations when the distribution of the unobserved vector is degenerate. In many applications this may be problematic as a natural interpretation for the unobserved process may be lost.

\subsection{The Structural EM algorithm}

Suppose that given a random sample of size $n$, we are interested in finding the tree that maximizes the likelihood function over all \emph{fully-observed} tree models. For every tree the maximum likelihood estimate is easily obtained; see~\cite[Section 4.4.2 and Section 5.2.1]{lauritzen:96}. In the naive approach we can search over all possible trees with a given number of vertices and find the one that gives the highest value of the likelihood function. The Chow--Liu algorithm~\cite{chow:1968} is a remarkably simple algorithm, which gives an efficient way to find the best tree approximation that maximizes the likelihood. 

The original approach was proposed for discrete data but it can be easily extended to the Gaussian case~\cite{tan:2011}. It only uses the fact that the MLE decomposes according to a tree so it can be also used in other similar scenarios to find best BIC and AIC trees~\cite{edwards2010selecting}; see~\cite{hojsgaard2012graphical} for further references and implementations. This approach boils down to using the Kruskal's algorithm to find the maximum cost spanning tree of a complete graph with edge weights given by sample mutual informations \cite{kruskal1956shortest}:
\setlength{\leftmargini}{2em}
\begin{itemize}
\item[1.] Compute all mutual informations of the sample distribution and order them from the largest to the smallest.
\item[2.] Move along this ordered sequence adding subsequently the corresponding edges unless adding an edge introduces a cycle. Stop when no more edges can be added.
\end{itemize}
If all mutual informations are different, then there is a unique best solution. If some of the weights are equal, then multiple solutions are possible, but they will all give the same value of the likelihood function. In the Gaussian case the mutual information is a simple monotone function of the corresponding correlations squared. Therefore, the same tree will be obtained after replacing mutual informations in step 1 above with squares of sample correlations. 

A natural idea is to use the Chow--Liu algorithm to learn latent tree models by using an EM-type algorithm. The \emph{structural EM algorithm}~\cite{friedman2002structural} is a numerical procedure to maximize the likelihood simultaneously over the continuous parameter $\theta$ and the discrete parameter $T$. It starts from a given parameter value and moves at each step strictly increasing the likelihood unless it is already in a local optimum. The E-step of the algorithm is the standard E-step of the EM algorithm. The M-step follows essentially the Chow--Liu algorithm with several modifications. A minor problem with the Chow--Liu algorithm, when used in the EM algorithm, is that it does not get any information about which vertices represent unobserved variables and so it often outputs a tree whose leaves can be potentially unobserved or internal vertices that can represent observed random variables. In this case it is easy to provide a tree, with leaves precisely corresponding to observed vertices, that gives the same observed likelihood. This procedure is described in more detail for the discrete case in~\cite[Section 5]{friedman2002structural}. We describe the general idea in the Gaussian case. Suppose that the Chow--Liu algorithm outputs a tree $T'$ with edge correlations $\rho_{uv}$, then:
\begin{itemize}
\item Remove all degree one vertices that represent unobserved vertices.
	\item If there is an induced chain $\overset{i_1}{*}-\overset{i_2}{\circ}-\overset{i_3}{\circ}-\cdots-\overset{i_{k-1}}{\circ}-\overset{i_k}{*}$, where $*$ stands for any vertex of degree at least three, then replace this chain with a single edge between $i_1$ and $i_k$. Set the correlation $\rho_{i_1 i_k}$ equal to the product $\rho_{i_1 i_2}\cdots \rho_{i_{k-1} i_k}$.
	\item If there is an internal vertex $v$ representing an observed random variable then add an auxiliary copy $v'$ of $v$ and an edge $(v,v')$. Set $\rho_{vv'}=1$.
\end{itemize}
This operation gives a leaf-labeled tree that leads to the same observed likelihood as $T'$ and we take it as the output of the M-step. If we are interested only in trivalent trees, the above procedure can be easily modified so that the output is a trivalent tree; c.f. Theorem~\ref{th:binaryred}.

\subsection{Phylogenetic invariants}

Given latent tree model over $T$, we associate to it the set $\cI_{T}$ of all polynomials vanishing on $M(T,\cY)$. The polynomials are expressed in terms of correlations in the Gaussian case and in terms of the raw probabilities in the discrete case. Every polynomial $f\in \cI_{T}$ is called a \textit{phylogenetic invariant}. Equation (\ref{eq:quarteteqs}) provides an example of an equation that must hold for a Gaussian latent tree model over a quartet tree. It is a basic result in algebraic geometry that $\cI_T$ is finitely generated, that is, it admits a finite basis of polynomials $\{f_1,\ldots,f_r\}$ such that every polynomial in $\cI_T$ is a polynomial combination of the $f_i$. 

The basic idea behind application of phylogenetic invariants is as follows. Given $n$ independent observations of $X$ we compute the sample distribution $\hat{p}$, which, by the law of large numbers, converges almost surely to the true data generating distribution $p^*$. Because $f(p^{*})=0$ for every $f\in \cI_{T}$, for large $n$ also $f(\hat{p})\approx 0$. The methods proposed in the phylogenetic  literature are mainly simple diagnostic tests that work with a given fixed finite set of invariants in $\cI_T$, which do not necessarily generate $\cI_T$. There is now considerable literature on the method of phylogenetic invariants and for many models all defining polynomials are understood. We refer to~\cite{allman_invariants2007,sumner2016developing} for an overview. 

The advantage of this approach is that the method based on phylogenetic invariants does not require parameter estimation. The disadvantage is that to large extent statistical theory behind their use is lacking. The method of phylogenetic invariants gives a way to select the best tree under a given criterion. It does not, however, give a way of quantifying how well the chosen tree fits the data because, in general, the distribution of phylogenetic invariants is too hard to analyze. 

Recently there has been some effort aimed at organizing the statistical theory behind phylogenetic invariants. For example, it has been observed by many authors that not all invariants are equally important and from the statistical point of view there is no sense in working with a full generating set of $\cI_T$. For example, invariants linking directly to tree metrics, called the \emph{edge invariants}, tend to be more robust with respect to the sampling error and are enough to distinguish between different models~\cite{casanellas2011relevant}. To test the edge invariants, \cite{eriksson2007uip} uses the singular value decomposition and the Frobenius norm to compute the distance of a matrix to the set of matrices of certain rank. Recently this method has been further improved by~\cite{fernandez2014invariant}. In its current form the method is robust and simulations show that it outperforms most of the commonly used methods.  

Focusing only on quadratic invariants allows us to generalise directly asymptotic chi-square tests for independence in a contingency table; see~\cite{sankoff1990designer}. In the Gaussian setting we can proceed as follows. For any four leaves of $T$, $ij/kl$ \emph{forms a quartet} in $T$ if the paths $\overline{ij}$ and $\overline{kl}$ have no vertices in common. It is clear from~(\ref{eq:RhoPath}) and~(\ref{eq:qij}) that for any quartet $ij/kl$ we have that $\rho_{ik}\rho_{jl}-\rho_{il}\rho_{jk}=0$ in the Gaussian case and $\tau_{ik}\tau_{jl}-\tau_{il}\tau_{jk} =0$ in the discrete case. For example~(\ref{eq:quarteteqs}) gives the equation that holds for the quartet tree in Example~\ref{ex:Gquartet}. Equations of this form are called tetrads~\cite{silva2006learning}. In the context of latent tree models using tetrads was suggested already by Judea Pearl~\cite[Section 8.3.5]{pearlBNbook} but with no statistical guidance. An example of a statistically guided quartet-based analysis was given recently for Gaussian latent tree models~\cite{ASSZ2014}.

Finally, the geometric description of latent tree models involves not only polynomial equalities but also polynomial inequalities. These inequalities cut out a large portion of the space described solely by equalities and therefore they should not be neglected~\cite{pwz-2009-semialgebraictrees}. The inequalities are harder to study than equalities but they are also understood well for general Markov models~\cite{allman2012semialgebraic,settimi2000gma,pwz-2009-semialgebraictrees}. Recently~\cite{ARSZ2013} showed how basic inequalities can be easily tested within the Bayesian framework to obtain a preliminary assessment of whether the data come from a Gaussian latent tree model.

\section{Discussion}

We gave a concise overview of the theory of latent tree models. We argued that, from the theoretical and the practical point of view, the general Markov models, and more generally, linear latent tree models, form the most important subfamily. We showed that for linear latent tree models the link to tree metrics provides a wide variety of learning procedures. The geometric viewpoint gives important insights but was not covered here in much detail. We refer to~\cite{LTbook} for further details. A more algorithmic overview of the latent tree model class is provided in~\cite{mourad2013survey}. We also skipped other research directions that we believe will become increasingly important in the coming years. Tensor representations of discrete latent tree models help to design learning procedures for latent tree models and more general graphical models with unobserved variables~\cite{anandkumar2014tensor,ishteva2013unfolding,ishteva2013hierarchical,parikh2011spectral}. In the numerical analysis community a closely related concept of hierarchical tensors has become popular \cite{hackbusch2009new}. There are also several recent approaches that introduce the nonparametric setting for latent tree models~\cite{song2014nonparametric,song2011kernel}.

\section*{Acknowledgements}

I would like to thank Anima Anandkumar, Sebastien Roch, and Nevin L. Zhang for helpful comments, clarifications, and literature suggestions. Comments from the anonymous referees allowed to substantially improve the manuscript.

\bibliographystyle{plain}
\bibliography{algebraic_statistics}

\begin{thebibliography}{100}

\bibitem{allman2009identifiability}
Elizabeth~S. Allman, Catherine Matias, and John~A. Rhodes.
\newblock Identifiability of parameters in latent structure models with many
  observed variables.
\newblock {\em Ann. Statist.}, 37(6A):3099--3132, 2009.

\bibitem{allman_invariants2007}
Elizabeth~S. Allman and John~A. Rhodes.
\newblock Phylogenetic invariants.
\newblock In {\em Reconstructing {E}volution: {N}ew {M}athematical and
  {C}omputational {A}dvances}, pages 108--146. Oxford Univ. Press, Oxford,
  2007.

\bibitem{allman2008pia}
Elizabeth~S. Allman and John~A. Rhodes.
\newblock Phylogenetic ideals and varieties for the general {M}arkov model.
\newblock {\em Adv. in Appl. Math.}, 40(2):127--148, 2008.

\bibitem{ARSZ2013}
Elizabeth~S. Allman, John~A. Rhodes, Bernd Sturmfels, and Piotr Zwiernik.
\newblock Tensors of nonnegative rank two.
\newblock {\em Linear Algebra Appl.}, 473:37--53, 2015.

\bibitem{allman2012semialgebraic}
Elizabeth~S. Allman, John~A. Rhodes, and Amelia Taylor.
\newblock A semialgebraic description of the general {M}arkov model on
  phylogenetic trees.
\newblock {\em SIAM J. Discrete Math.}, 28(2):736--755, 2014.

\bibitem{anandkumar2011spectral}
Animashree Anandkumar, Kamalika Chaudhuri, Daniel~J Hsu, Sham~M Kakade,
  Le~Song, and Tong Zhang.
\newblock Spectral methods for learning multivariate latent tree structure.
\newblock In {\em Advances in Neural Information Processing Systems}, pages
  2025--2033, 2011.

\bibitem{anandkumar2014tensor}
Animashree Anandkumar, Rong Ge, Daniel Hsu, Sham~M. Kakade, and Matus
  Telgarsky.
\newblock Tensor decompositions for learning latent variable models.
\newblock {\em J. Mach. Learn. Res.}, 15:2773--2832, 2014.

\bibitem{anderson1956statistical}
T.~W. Anderson and Herman Rubin.
\newblock Statistical inference in factor analysis.
\newblock In {\em Proceedings of the {T}hird {B}erkeley {S}ymposium on
  {M}athematical {S}tatistics and {P}robability, 1954--1955, vol. {V}}, pages
  111--150. University of California Press, Berkeley and Los Angeles, 1956.

\bibitem{barry1987statistical}
Daniel Barry and J.~A. Hartigan.
\newblock Statistical analysis of hominoid molecular evolution.
\newblock {\em Statist. Sci.}, 2(2):191--210, 1987.
\newblock With comments by Stephen Portnoy and Joseph Felsenstein and a reply
  by the authors.

\bibitem{bishop2006pattern}
Christopher~M Bishop.
\newblock {\em Pattern recognition and machine learning}.
\newblock Springer, 2006.

\bibitem{bishop1998hierarchical}
Christopher~M Bishop and Michael~E Tipping.
\newblock A hierarchical latent variable model for data visualization.
\newblock {\em IEEE Transactions on Pattern Analysis and Machine Intelligence},
  20(3):281--293, 1998.

\bibitem{buneman1971}
Peter Buneman.
\newblock The recovery of trees from measures of dissimilarity.
\newblock In F.~Hodson et~al., editor, {\em Mathematics in the Archaeological
  and Historical Sciences}, pages 387--395. Edinburgh University Press, 1971.

\bibitem{casanellas2011relevant}
Marta Casanellas and Jes{\'u}s Fern{\'a}ndez-S{\'a}nchez.
\newblock Relevant phylogenetic invariants of evolutionary models.
\newblock {\em J. Math. Pures Appl. (9)}, 96(3):207--229, 2011.

\bibitem{castro2004network}
Rui Castro, Mark Coates, Gang Liang, Robert Nowak, and Bin Yu.
\newblock Network tomography: recent developments.
\newblock {\em Statistical Science}, pages 499--517, 2004.

\bibitem{cavallisforza1967pam}
LL~Cavalli-Sforza and AWF Edwards.
\newblock {Phylogenetic analysis: models and estimation procedures}.
\newblock {\em Evolution}, 21(3):550--570, 1967.

\bibitem{chang1996frm}
Joseph~T. Chang.
\newblock {Full reconstruction of {M}arkov models on evolutionary trees:
  Identifiability and consistency}.
\newblock {\em Mathematical Biosciences}, 137(1):51--73, 1996.

\bibitem{chang1991reconstruction}
Joseph~T Chang and John~A Hartigan.
\newblock Reconstruction of evolutionary trees from pairwise distributions on
  current species.
\newblock In {\em Computing science and statistics: Proceedings of the 23rd
  symposium on the interface}, volume 254, page 257. Interface Foundation,
  Fairfax Station, VA, 1991.

\bibitem{chen2015progressive}
Peixian Chen, Nevin~L. Zhang, Leonard Poon, and Zhourong Chen.
\newblock Progressive {EM} for latent tree models and hierarchical topic
  detection.
\newblock In {\em AAAI}, pages 1498--1504, 2016.

\bibitem{chen2012model}
Tao Chen, Nevin~L Zhang, Tengfei Liu, Kin~Man Poon, and Yi~Wang.
\newblock Model-based multidimensional clustering of categorical data.
\newblock {\em Artificial Intelligence}, 176(1):2246--2269, 2012.

\bibitem{choi2001multiscale}
Hyeokho Choi and Richard~G. Baraniuk.
\newblock Multiscale image segmentation using wavelet-domain hidden {M}arkov
  models.
\newblock {\em IEEE Trans. Image Process.}, 10(9):1309--1321, 2001.

\bibitem{choi2010exploiting}
Myung~Jin Choi, Joseph~J Lim, Antonio Torralba, and Alan~S. Willsky.
\newblock Exploiting hierarchical context on a large database of object
  categories.
\newblock In {\em Computer Vision and Pattern Recognition (CVPR)}, pages
  129--136. IEEE, 2010.

\bibitem{choi:2011}
Myung~Jin Choi, Vincent Y.~F. Tan, Animashree Anandkumar, and Alan~S. Willsky.
\newblock Learning latent tree graphical models.
\newblock {\em J. Mach. Learn. Res.}, 12:1771--1812, 2011.

\bibitem{chor2000mml}
Benny Chor, Michael~D. Hendy, Barbara~R. Holland, and David Penny.
\newblock Multiple maxima of likelihood in phylogenetic trees: An analytic
  approach.
\newblock {\em Molecular Biology and Evolution}, 17(10):1529--1541, 2000.

\bibitem{chow:1968}
C.~K. Chow and C.~N. Liu.
\newblock Approximating discrete probability distributions with dependence
  trees.
\newblock {\em IEEE Trans. Inform. Theory}, 14:462--467, 1968.

\bibitem{crouse1998wavelet}
Matthew~S. Crouse, Robert~D. Nowak, and Richard~G. Baraniuk.
\newblock Wavelet-based statistical signal processing using hidden {M}arkov
  models.
\newblock {\em IEEE Trans. Signal Process.}, 46(4):886--902, 1998.

\bibitem{daskalakis2011phylogenies}
Constantinos Daskalakis, Elchanan Mossel, and Sebastien Roch.
\newblock Phylogenies without branch bounds: Contracting the short, pruning the
  deep.
\newblock {\em SIAM Journal on Discrete Mathematics}, 25(2):872--893, 2011.

\bibitem{DLWZ2014}
Mathias Drton, Shaowei Lin, Luca Weihs, and Piotr Zwiernik.
\newblock Marginal likelihood and model selection for {G}aussian latent tree
  and forest models.
\newblock {\em Bernoulli}, 23(2):1202--1232, 2017.

\bibitem{drton2013sBIC}
Mathias Drton and Martyn Plummer.
\newblock A {B}ayesian information criterion for singular models.
\newblock {\em Journal of the Royal Statistical Society: Series B (Statistical
  Methodology)}, 79(2):323--380, 2017.

\bibitem{oberwolfach2009}
Mathias Drton, Bernd Sturmfels, and Seth Sullivant.
\newblock {\em Lectures on Algebraic Statistics}.
\newblock Oberwolfach Seminars Series. Birkhauser Verlag AG, 2009.

\bibitem{durbin1998bsa}
Richard Durbin, Anders Krogh, Graeme Mitchison, and Sean~R. Eddy.
\newblock {\em {Biological Sequence Analysis: Probabilistic Models of Proteins
  and Nucleic Acids}}.
\newblock Cambridge University Press, 1998.

\bibitem{edwards2010selecting}
David Edwards, Gabriel~CG De~Abreu, and Rodrigo Labouriau.
\newblock Selecting high-dimensional mixed graphical models using minimal aic
  or bic forests.
\newblock {\em BMC bioinformatics}, 11(1):1, 2010.

\bibitem{erdos1999few}
P{\'e}ter~L. Erd{\H{o}}s, Michael~A. Steel, L{\'a}szl{\'o}~A. Sz{\'e}kely, and
  Tandy~J. Warnow.
\newblock A few logs suffice to build (almost) all trees (part 1).
\newblock {\em Random Structures and Algorithms}, 14(2):153--184, 1999.

\bibitem{erdos199977}
P{\'e}ter~L. Erd{\H{o}}s, Michael~A. Steel, L{\'a}szl{\'o}~A. Sz{\'e}kely, and
  Tandy~J. Warnow.
\newblock A few logs suffice to build (almost) all trees (part 2).
\newblock {\em Theoretical Computer Science}, 221(1):77 -- 118, 1999.

\bibitem{eriksson2010toward}
Brian Eriksson, Gautam Dasarathy, Paul Barford, and Robert Nowak.
\newblock Toward the practical use of network tomography for internet topology
  discovery.
\newblock In {\em INFOCOM, 2010 Proceedings IEEE}, pages 1--9. IEEE, 2010.

\bibitem{eriksson2007uip}
Nicholas Eriksson.
\newblock {\em {Using Invariants for Phylogenetic Tree Construction}}, volume
  149 of {\em The IMA Volumes in Mathematics and Its Applications}, pages
  89--108.
\newblock Springer, 2007.

\bibitem{evans2015graphs}
Robin~J. Evans.
\newblock Graphs for margins of {B}ayesian networks.
\newblock {\em Scandinavian Journal of Statistics}, 43(3):625--648, 2016.
\newblock 10.1111/sjos.12194.

\bibitem{felsenstein2004inferring}
Joseph Felsenstein.
\newblock {\em Inferring phylogenies}.
\newblock Sinauer Associates Sunderland, 2004.

\bibitem{fernandez2014invariant}
Jes{\'u}s Fern{\'a}ndez-S{\'a}nchez and Marta Casanellas.
\newblock Invariant versus classical approaches when evolution is heterogeneous
  across sites and lineages.
\newblock {\em arXiv:1405.6546}, 2014.

\bibitem{friedman2002structural}
Nir Friedman, Matan Ninio, Itsik Pe'er, and Tal Pupko.
\newblock A structural {EM} algorithm for phylogenetic inference.
\newblock {\em Journal of Computational Biology}, 9(2):331--353, 2002.

\bibitem{frydenberg1990}
Morten Frydenberg.
\newblock The chain graph {M}arkov property.
\newblock {\em Scand. J. Statist.}, 17(4):333--353, 1990.

\bibitem{guardiano2005parametric}
Cristina Guardiano and Giuseppe Longobardi.
\newblock Parametric comparison and language taxonomy.
\newblock In Montserrat Batllori, Maria-Llu{\"\i}sa Hernanz, Carme Picallo, and
  Francesc Roca, editors, {\em Grammaticalization and Parametric Variation},
  pages 149--174. Oxford University Press, Aug 2005.

\bibitem{hackbusch2009new}
Wolfgang Hackbusch and Stefan K{\"u}hn.
\newblock A new scheme for the tensor representation.
\newblock {\em Journal of Fourier analysis and applications}, 15(5):706--722,
  2009.

\bibitem{hiscott2016efficient}
Gordon Hiscott, Colin Fox, Matthew Parry, and David Bryant.
\newblock Efficient recycled algorithms for quantitative trait models on
  phylogenies.
\newblock {\em Genome biology and evolution}, 8(5):1338--1350, 2016.

\bibitem{hojsgaard2012graphical}
S{\o}ren H{\o}jsgaard, David Edwards, and Steffen Lauritzen.
\newblock {\em Graphical models with R}.
\newblock Springer Science \& Business Media, 2012.

\bibitem{holland01012013}
Barbara~R. Holland, Peter~D. Jarvis, and Jeremy~G. Sumner.
\newblock Low-parameter phylogenetic inference under the general {M}arkov
  model.
\newblock {\em Systematic Biology}, 62(1):78--92, 2013.

\bibitem{huang2014scalable}
Furong Huang, U.N. Niranjan, Ioakeim Perros, Robert Chen, Jimeng Sun, and Anima
  Anandkumar.
\newblock Scalable latent tree model and its application to health analytics.
\newblock {\em arXiv:1406.4566}, 2014.

\bibitem{ishteva2013unfolding}
Mariya Ishteva, Haesun Park, and Le~Song.
\newblock Unfolding latent tree structures using 4th order tensors.
\newblock In {\em ICML (3)}, pages 316--324, 2013.

\bibitem{ishteva2013hierarchical}
Mariya Ishteva, L~Song, H~Park, A~Parikh, and E~Xing.
\newblock Hierarchical tensor decomposition of latent tree graphical models.
\newblock In {\em The 30th International Conference on Machine Learning (ICML
  2013)}, 2013.

\bibitem{jayaswal2007estimation}
Vivek Jayaswal, John Robinson, and Lars~S. Jermiin.
\newblock Estimation of phylogeny and invariant sites under the general
  {M}arkov model of nucleotide sequence evolution.
\newblock {\em Systematic Biology}, 56(2):155--162, 2007.

\bibitem{jukes1969epm}
Thomas~H. Jukes and Charles~R. Cantor.
\newblock {Evolution of protein molecules}.
\newblock {\em Mammalian Protein Metabolism}, 3:21--132, 1969.

\bibitem{kruskal1956shortest}
Joseph~B Kruskal.
\newblock On the shortest spanning subtree of a graph and the traveling
  salesman problem.
\newblock {\em Proceedings of the American Mathematical Society}, 7(1):48--50,
  1956.

\bibitem{lake1994ret}
James~A. Lake.
\newblock {Reconstructing evolutionary trees from DNA and protein sequences:
  Paralinear distances}.
\newblock {\em Proceedings of the National Academy of Sciences of the United
  States of America}, 91(4):1455--1459, 1994.

\bibitem{lauritzen:96}
Steffen~L. Lauritzen.
\newblock {\em Graphical {M}odels}, volume~17 of {\em Oxford Statistical
  Science Series}.
\newblock Oxford University Press, 1996.
\newblock Oxford Science Publications.

\bibitem{lawrence2004gaussian}
Neil~D. Lawrence.
\newblock Gaussian process latent variable models for visualisation of high
  dimensional data.
\newblock {\em Advances in Neural Information Processing Systems},
  16(3):329--336, 2004.

\bibitem{shaowei_rlct}
Shaowei Lin.
\newblock Ideal-theoretic strategies for asymptotic approximation of marginal
  likelihood integrals.
\newblock {\em Journal of Algebraic Statistics}, 8(1), 2017.

\bibitem{lockhart1994recovering}
Peter~J. Lockhart, Michael~A. Steel, Michael~D. Hendy, and David Penny.
\newblock Recovering evolutionary trees under a more realistic model of
  sequence evolution.
\newblock {\em Molecular Biology and Evolution}, 11(4):605--612, 1994.

\bibitem{makhijani2012accelerated}
Mahender~K Makhijani, Niranjan Balu, Kiyofumi Yamada, Chun Yuan, and Krishna~S
  Nayak.
\newblock Accelerated 3d merge carotid imaging using compressed sensing with a
  hidden {M}arkov tree model.
\newblock {\em Journal of Magnetic Resonance Imaging}, 36(5):1194--1202, 2012.

\bibitem{mossel2003impossibility}
Elchanan Mossel.
\newblock On the impossibility of reconstructing ancestral data and
  phylogenies.
\newblock {\em Journal of Computational Biology}, 10(5):669--676, 2003.

\bibitem{mossel2004phase}
Elchanan Mossel.
\newblock Phase transitions in phylogeny.
\newblock {\em Transactions of the American Mathematical Society},
  356(6):2379--2404, 2004.

\bibitem{mossel:2013}
Elchanan Mossel, S\'ebastien Roch, and Allan Sly.
\newblock Robust estimation of latent tree graphical models: Inferring hidden
  states with inexact parameters.
\newblock {\em IEEE Trans. Inform. Theory}, 59(7):4357--4373, 2013.

\bibitem{mourad2013survey}
Rapha{\"e}l Mourad, Christine Sinoquet, Nevin~Lianwen Zhang, Tengfei Liu,
  Philippe Leray, et~al.
\newblock A survey on latent tree models and applications.
\newblock {\em J. Artif. Intell. Res.(JAIR)}, 47:157--203, 2013.

\bibitem{paradis2011analysis}
Emmanuel Paradis.
\newblock {\em Analysis of Phylogenetics and Evolution with R}.
\newblock Springer Science \& Business Media, 2011.

\bibitem{parikh2011spectral}
Ankur~P Parikh, Le~Song, and Eric~P Xing.
\newblock A spectral algorithm for latent tree graphical models.
\newblock In {\em Proceedings of the 28th International Conference on Machine
  Learning (ICML-11)}, pages 1065--1072, 2011.

\bibitem{pearlBNbook}
Judea Pearl.
\newblock {\em {Probabilistic Reasoning in Intelligent Systems: Networks of
  Plausible Inference}}.
\newblock The Morgan Kaufmann Series in Representation and Reasoning. Morgan
  Kaufmann, San Mateo, CA, 1988.

\bibitem{pearl2000}
Judea Pearl.
\newblock {\em Causality: {M}odels, {R}easoning, and {I}nference}.
\newblock Cambridge University Press, New York, NY, USA, 2000.

\bibitem{pfeiffer:CNE23820}
Michael Pfeiffer, Marion Betizeau, Julie Waltispurger, Sabina~Sara Pfister,
  Rodney~J. Douglas, Henry Kennedy, and Colette Dehay.
\newblock Unsupervised lineage-based characterization of primate precursors
  reveals high proliferative and morphological diversity in the {OSVZ}.
\newblock {\em Journal of Comparative Neurology}, 524(3):535--563, 2016.

\bibitem{rabiner1989tutorial}
Lawrence~R Rabiner.
\newblock A tutorial on hidden {M}arkov models and selected applications in
  speech recognition.
\newblock {\em Proceedings of the IEEE}, 77(2):257--286, 1989.

\bibitem{ringe2002indo}
Don Ringe, Tandy Warnow, and Ann Taylor.
\newblock Indo-european and computational cladistics.
\newblock {\em Transactions of the Philological Society}, 100(1):59--129, 2002.

\bibitem{roch2006short}
Sebastien Roch.
\newblock A short proof that phylogenetic tree reconstruction by maximum
  likelihood is hard.
\newblock {\em IEEE/ACM Trans. Comput. Biol. Bioinformatics}, 3(1):92--,
  January 2006.

\bibitem{roch2015phase}
Sebastien Roch and Allan Sly.
\newblock Phase transition in the sample complexity of likelihood-based
  phylogeny inference.
\newblock {\em arXiv:1508.01964}, 2015.

\bibitem{romberg2001bayesian}
Justin~K Romberg, Hyeokho Choi, and Richard~G Baraniuk.
\newblock Bayesian tree-structured image modeling using wavelet-domain hidden
  {M}arkov models.
\newblock {\em IEEE Transactions on image processing}, 10(7):1056--1068, 2001.

\bibitem{sankoff1990designer}
David Sankoff.
\newblock {Designer invariants for large phylogenies}.
\newblock {\em Molecular Biology and Evolution}, 7(3):255, 1990.

\bibitem{schwarz1978edm}
Gideon Schwarz.
\newblock {Estimating the dimension of a model}.
\newblock {\em Annals of Statistics}, 6(2):461--464, 1978.

\bibitem{semple2003pol}
Charles Semple and Mike Steel.
\newblock {\em Phylogenetics}, volume~24 of {\em Oxford Lecture Series in
  Mathematics and Its Applications}.
\newblock Oxford University Press, Oxford, 2003.

\bibitem{settimi2000gma}
Raffaella Settimi and Jim~Q. Smith.
\newblock Geometry, moments and conditional independence trees with hidden
  variables.
\newblock {\em Ann. Statist.}, 28(4):1179--1205, 2000.

\bibitem{shiers2014gaussian}
Nathaniel Shiers, John~AD Aston, Jim~Q Smith, and John~S Coleman.
\newblock Gaussian tree constraints applied to acoustic linguistic functional
  data.
\newblock {\em Journal of Multivariate Analysis}, 154:199--215, 2017.

\bibitem{ASSZ2014}
Nathaniel Shiers, Piotr Zwiernik, John~A. Aston, and James~Q. Smith.
\newblock The correlation space of {G}aussian latent tree models and model
  selection without fitting.
\newblock {\em Biometrika}, 2016.

\bibitem{shu2016syntactic}
Kevin Shu, Sharjeel Aziz, Vy-Luan Huynh, David Warrick, and Matilde Marcolli.
\newblock Syntactic phylogenetic trees.
\newblock {\em arXiv:1607.02791}, 2016.

\bibitem{silva2006learning}
Ricardo Silva, Richard Scheine, Clark Glymour, and Peter Spirtes.
\newblock Learning the structure of linear latent variable models.
\newblock {\em Journal of Machine Learning Research}, 7(Feb):191--246, 2006.

\bibitem{song2014nonparametric}
Le~Song, Animashree Anandkumar, Bo~Dai, and Bo~Xie.
\newblock Nonparametric estimation of multi-view latent variable models.
\newblock In {\em Proceedings of the 31st International Conference on Machine
  Learning (ICML-14)}, pages 640--648, 2014.

\bibitem{song2011kernel}
Le~Song, Eric~P. Xing, and Ankur~P. Parikh.
\newblock Kernel embeddings of latent tree graphical models.
\newblock In J.~Shawe-Taylor, R.~S. Zemel, P.~L. Bartlett, F.~Pereira, and
  K.~Q. Weinberger, editors, {\em Advances in Neural Information Processing
  Systems 24}, pages 2708--2716. Curran Associates, Inc., 2011.

\bibitem{steel2004}
M.~Steel.
\newblock Recovering a tree from the leaf colourations it generates under a
  {M}arkov model.
\newblock {\em Applied Mathematics Letters}, 7(2):19 -- 23, 1994.

\bibitem{steel2001my}
M~Steel.
\newblock My favourite conjecture.
\newblock {\em Preprint}, 2001.

\bibitem{steel1993invertible}
M.~Steel, M.D. Hendy, and D.~Penny.
\newblock Invertible models of sequence evolution.
\newblock Mathematical and Information Science report 93/02, Massey University,
  1993.

\bibitem{sullivant2008agg}
Seth Sullivant.
\newblock {Algebraic geometry of Gaussian Bayesian networks}.
\newblock {\em Advances in Applied Mathematics}, 40(4):482--513, 2008.

\bibitem{sumner2016developing}
Jeremy~G Sumner, Amelia Taylor, Barbara~R Holland, and Peter~D Jarvis.
\newblock Developing a statistically powerful measure for quartet tree
  inference using phylogenetic identities and {M}arkov invariants.
\newblock {\em Journal of Mathematical Biology}, pages 1--36, 2017.

\bibitem{tan:2011}
Vincent Y.~F. Tan, Animashree Anandkumar, and Alan~S. Willsky.
\newblock Learning high-dimensional {M}arkov forest distributions: analysis of
  error rates.
\newblock {\em J. Mach. Learn. Res.}, 12:1617--1653, 2011.

\bibitem{thurstone1934vectors}
Louis~L. Thurstone.
\newblock The vectors of mind.
\newblock {\em Psychological Review}, 41(1):1, 1934.

\bibitem{vermapearl91}
Thomas~S. Verma and Judea Pearl.
\newblock Equivalence and {S}ynthesis of {C}ausal {M}odels.
\newblock In Piero~P. Bonissone, Max Henrion, Laveen~N. Kanal, and John~F.
  Lemmer, editors, {\em UAI '90: Proceedings of the Sixth Annual Conference on
  Uncertainty in Artificial Intelligence, MIT, Cambridge, MA, USA, July 27-29,
  1990}. Elsevier, October 1991.

\bibitem{wainwright2008graphical}
Martin~J. Wainwright and Michael~I. Jordan.
\newblock Graphical models, exponential families, and variational inference.
\newblock {\em Foundations and Trends in Machine Learning}, 1(1-2):1--305,
  2008.

\bibitem{wang2006severity}
Yi~Wang and Nevin~L. Zhang.
\newblock Severity of local maxima for the {EM} algorithm: Experiences with
  hierarchical latent class models.
\newblock In {\em Probabilistic Graphical Models}, pages 301--308. Citeseer,
  2006.

\bibitem{watanabe_book}
Sumio Watanabe.
\newblock {\em Algebraic Geometry and Statistical Learning Theory}.
\newblock Number~25 in Cambridge Monographs on Applied and Computational
  Mathematics. Cambridge University Press, 2009.
\newblock ISBN-13: 9780521864671.

\bibitem{watanabe2013widely}
Sumio Watanabe.
\newblock A widely applicable {B}ayesian information criterion.
\newblock {\em Journal of Machine Learning Research}, 14(Mar):867--897, 2013.

\bibitem{willsky2002multiresolution}
Alan~S. Willsky.
\newblock Multiresolution {M}arkov models for signal and image processing.
\newblock {\em Proceedings of the IEEE}, 90(8):1396--1458, 2002.

\bibitem{vzabokrtsky2009hidden}
Zden{\v{e}}k {\v{Z}}abokrtsk{\`y} and Martin Popel.
\newblock Hidden {M}arkov tree model in dependency-based machine translation.
\newblock In {\em Proceedings of the ACL-IJCNLP 2009 Conference Short Papers},
  pages 145--148. Association for Computational Linguistics, 2009.

\bibitem{zhangCluster}
Nevin~L. Zhang.
\newblock Hierarchical latent class models for cluster analysis.
\newblock {\em J. Mach. Learn. Res.}, 5:697--723, 2003/04.

\bibitem{zhang2017latent}
Nevin~L Zhang and Leonard~KM Poon.
\newblock Latent tree analysis.
\newblock In {\em AAAI}, pages 4891--4898, 2017.

\bibitem{zou06042011}
Liwen Zou, Edward Susko, Chris Field, and Andrew~J. Roger.
\newblock The parameters of the {B}arry and {H}artigan general {M}arkov model
  are statistically nonidentifiable.
\newblock {\em Systematic Biology}, 2011.

\bibitem{pwz-2010-bic}
Piotr Zwiernik.
\newblock Asymptotic behaviour of the marginal likelihood for general {M}arkov
  models.
\newblock {\em J. Mach. Learn. Res.}, 12:3283--3310, 2011.

\bibitem{LTbook}
Piotr Zwiernik.
\newblock {\em Semialgebraic statistics and latent tree models}, volume 146 of
  {\em Monographs on Statistics and Applied Probability}.
\newblock Chapman \& Hall/CRC, 2016.

\bibitem{pwz-2009-semialgebraictrees}
Piotr Zwiernik and Jim~Q. Smith.
\newblock Implicit inequality constraints in a binary tree model.
\newblock {\em Electron. J. Statist.}, 5:1276--1312, 2011.

\bibitem{pwz-2010-identifiability}
Piotr Zwiernik and Jim~Q. Smith.
\newblock {Tree-cumulants and the geometry of binary tree models}.
\newblock {\em Bernoulli}, 18(1):290--321, January 2012.

\end{thebibliography}
\bigskip
\textbf{Author's address:}\\[.2cm] Piotr Zwiernik, Universitat Pompeu Fabra, Department of Economics and Business, Ramon Trias Fargas 25-27, 
08005 Barcelona, Spain.\\ E-mail: \texttt{piotr.zwiernik@upf.edu}

\end{document}